\documentclass[11pt, a4paper]{article}
\usepackage{amsmath,amssymb,amsfonts}
\usepackage{a4wide}
\usepackage{psfig,epsfig,epic,eepic,graphicx}
\newtheorem{theorem}{Theorem}
\newtheorem{lemma}[theorem]{Lemma}
\newtheorem{proposition}[theorem]{Proposition}

\makeatletter

\@addtoreset{equation}{section}
\makeatother
\newcommand{\E}{{\mathbb E}}
\newcommand{\T}{{\mathbb T}}
\newcommand{\N}{{\mathbb N}}

\newcommand{\Q}{{\mathbb Q}}
\newcommand{\R}{{\mathbb R}}

\newcommand{\pa}{{\partial}}
\newcommand{\na}{{\nabla}}
\newcommand{\dis}{{\displaystyle}}
\newcommand{\eps}{{\varepsilon}}

\def\div{\hbox{div  }}

\def\dis{\displaystyle}
\title{Wall laws for fluid flows   at a \\   boundary with random roughness}
\author{\footnote{DMA, Ecole Normale Sup\'erieure,
 45 rue d'Ulm,75005 Paris} Arnaud {\sc Basson}, 
 \footnote{DMA/CNRS, Ecole    Normale Sup\'erieure, 45 rue d'Ulm,75005
   Paris}  David {\sc G\'erard-Varet}}

\begin{document}
\maketitle

The general concern of this paper is the effect of
  rough boundaries on fluids. We consider a stationary flow, governed
  by incompressible Navier-Stokes equations, in an infinite 
 domain bounded by two  horizontal rough plates. The roughness is modeled
 by a spatially homogeneous random field, with characteristic size
 $\eps$. A mathematical analysis of the flow for small $\eps$ is
 performed. The  Navier's wall law is rigorously deduced from this
 analysis. This extends substantially former results obtained in the case of
 periodic roughness, notably in \cite{Jager:2001,Jager:2003}. 

\tableofcontents
%%%%%%%%%%%%%%%%%%%%%%%%%%%%%%%%%%%%%%%%%%%%%%%%%%%%%%%%%%%%%%%%%%%%%%%%%%%
%%%%%%%%%%%%%%%%%%%%%%%%%%%%%%%%%%%%%%%%%%%%%%%%%%%%%%%%%%%%%%%%%%%%%%%%%%%
\section{Introduction}

The understanding of roughness-induced effects is
 a major concern in fluid dynamics. Indeed, many examples of physical
 relevance involve rough boun\-da\-ries. By ``rough'', we mean that the 
spatial  variations are small compared to the typical length of the problem. 
For instance, in geophysics, the bottom of the oceans and the shores are rough 
with respect to the large scale flow. Also, in an 
industrial framework, containers have often imperfections that qualify them
as rough. 

The main problem is to know in which way such irregular boundaries 
affect the flow. This is especially important with regards to numerical
computations: indeed,  roughness is in general  too small to be captured by the
discretization grid of the simulations. 

To overcome this difficulty, one often relies on {\em wall laws}. A wall law is
a boundary condition that is imposed on an artificial boundary  inside the
domain. The idea is to filter out the precise description of the flow near
the real rough boundary. The wall law should only reflect the large scale
effect of the roughness, in the spirit of a homogeneization process. 

In many cases, the determination of wall laws relies on formal calculations, 
grounded by empirical arguments (see for instance \cite{Bechert:1989,
  Luchini:1995}).  
{\em The present paper is a mathematical justification of some wall laws, 
in the  case of an incompressible viscous fluid}. We  consider
 the two-dimensional stationary Navier-Stokes equations: 
\begin{equation} \label{NS}
\left\{
\begin{aligned}
& u \cdot \na u + \na p - \nu  \Delta u = 0, \: x \in \Omega^\eps,\\
& \div u = 0, \: x \in \Omega^\eps,
\end{aligned}
\right.
\end{equation} 
in a domain $\Omega^\eps$ of channel type: 
\begin{equation*}  
\Omega^\eps  \: = \: \left\{ (x_1, x_2) \in \R^2, \:
 \gamma_l^{\eps}(x_1)  < x_2 < \gamma_u^{\eps}(x_1)
\right\}, 
\end{equation*}
where the lower and upper boundaries $\gamma_l^\eps$ and
$\gamma_u^\eps$ are to be precised. As usual, the fields $u = (u_1, u_2)(x)
  \in \R^2, \:$ $p = p(x) \in \R$
  are the velocity and the pressure, and $\nu >
  0$ is the kinematic viscosity. Equations
  \eqref{NS} are completed with the classical no-slip conditions  
\begin{equation} \label{BC} 
 u\vert_{\pa \Omega^\eps} = 0.
\end{equation} 
Moreover, we prescribe the fluid flux through the channel, that is
condition
\begin{equation} \label{flux}
\int_{\sigma(x_1)} \!\!\!\! u_1 \:= \: \phi,
\end{equation}
where $\phi > 0$ is a given constant flux, and
$$\sigma(x_1) = \left\{
\gamma_l^{\eps}(x_1)  < x_2 < \gamma_u^{\eps}(x_1) \right\}$$
is a vertical section of the channel at $x_1$.
 Remark that by  incompressibility and
boundary condition \eqref{BC}, the left hand-side of \eqref{flux} does not
depend on $x_1$. The functions  $\gamma_l^\eps$ and
$\gamma_u^\eps$  model rough plates, with small characteristic size
$\eps$. Broadly, they read: 
$$ \gamma_l^\eps = -\eps
\gamma_l(x_1/\eps), \quad  
\gamma_u^\eps = 1 + \eps 
\gamma_u(x_1/\eps),  $$
for Lipschitz functions
$$\gamma_l = \gamma_l(y_1) \in
]0,1[, \quad \gamma_u = \gamma_u(y_1) \in ]0,1[.$$  
More precise  assumptions on  $\gamma_l$ and $\gamma_u$ will be
made further on. The set $\Omega = \R \times ]0,1[$ will be called the
    {\em interior domain}. 

We wish to study solutions $u^\eps$ of \eqref{NS}, \eqref{BC},
\eqref{flux}, and to determine  appropriate wall laws for this  system.
 In other words,  we look for 
 operators  ${\cal B}^\eps(x, D_x)$ such that solutions $v^\eps$ of the
 {\em interior system}  
\begin{equation} \label{reduced}
     \left\{
\begin{aligned}
& v \cdot \na v + \na q - \nu  \Delta v = 0, \: x \in \Omega, \\
& \div v = 0, \: x \in \Omega, \\
& \int_{\sigma(x_1)} \!\!\!\! v_1 \:= \: \phi, \\
& {\cal B}^\eps(x, D_x)(v) \vert_{\pa \Omega} = 0, 
\end{aligned}
\right.
\end{equation} 
approximate well  $u^\eps$ in $\Omega$, for small $\eps$.

The  mathematical treatment  of wall laws has been the matter of many
articles. The note \cite{Achdou:1995} is devoted to the analysis of  Laplace
equation in an annular domain with perforations. Numerical and formal
computations for fluid flows can be found in \cite{Achdou:1998a}
\cite{Achdou:1998}. An analysis of  Couette  flows in rough domains 
 has been performed in \cite{Amirat:2001a}. Let us also mention the
 important contributions of J\"ager and Mikelic on wall laws for 
channel flows (\cite{Jager:2001,Jager:2003}).  We 
refer to \cite{Jager:2000}
on a related problem with porous boundaries. Finally, see 
\cite{Gerard-Varet:2003b}, \cite{Bresch:2005} for study of 
roughness-induced effects on some geophysical systems. 

All these articles are devoted to {\em periodic roughness},
meaning that the boundary functions $\gamma_l$ and $\gamma_u$ are
periodic. This is of course a mathematical simplification, which is highly
unrealistic from the point of view of physics. {\em Our goal here is to
  drop this restriction, and treat non-periodic roughness. Precisely, we 
consider roughness that is distributed following a spatially homogeneous
random field.} A complete description of the rough domain will be given in
the next section.  

Following \cite{Jager:2001} in the periodic case,  special attention is
 paid  to  the simple Dirichlet wall law:
\begin{equation} \label{Dirichlet}
 {\cal B}^\eps(x, D_x)(v) \vert_{\pa \Omega} = v \vert_{\pa \Omega} = 0, 
\end{equation}

and to the Navier's friction law:
\begin{equation} \label{Navier}
 {\cal B}^\eps(x, D_x)(v) \vert_{\pa \Omega} = \left( C_{\eps}(x) \frac{\pa
  v_\tau}{\pa n} - v_\tau \right)\vert_{\pa \Omega} = 0.  
\end{equation}
introduced by Navier \cite{Navier:1827} and extensively used in simulations
of geophysical flows.
Losely, we  show two main results:
\begin{enumerate}
\item  The Dirichlet wall law yields a
$O(\eps)$ approximation of the real solution $u^\eps$, that is $u^\eps -
  v^\eps$ is $O(\eps)$ in an appropriate quadratic norm, to be described in
  the next section. 
\item  For appropriate $C_\eps$, the Navier's law yields a $o(\eps)$
  approximation of the real solution $u^\eps$. 
\end{enumerate}
These results extend those of \cite{Jager:2001}. They are deduced for a
precise description of $u^\eps$ for small $\eps$, especially of the boundary
layer flow near $\pa \Omega^\eps$. Precise statements, including the 
expression of $C_\eps$, will be given in the next section. 

To end this introduction, let us point out some difficulties related to the
proof of these results.
 First, we consider a domain $\Omega^\eps$ that is not bounded in the
tangential direction ($x_1 \in \R$). To our knowledge, previous studies
dealt with bounded channel domains, wether with  lateral boundaries (plus
in- and out-flux lateral boundary conditions, see \cite{Jager:2001})
 or with periodic boundary conditions, see \cite{Amirat:2001a}. Note that
 such periodicity condition is not compatible with our non-periodic 
 roughness. Due to  the unbounded channel domain, we work with 
only  locally integrable functions, which leads to completely 
different treatment of
the energy estimates.  Secondly, as the roughness is  non-periodic, the 
boundary layer system is  more complex. Due to the lack of compactness both 
in the tangential  and transverse variables, we are not able to solve it in
 a deterministic  setting. We use a variational formulation that involves
 the random variable. In addition to this problem, the behaviour of the
 boundary layer profile far from the boundary is not obvious. This can be
 understood using  formally the tangential Fourier transform. Indeed,
 in the periodic setting, the Fourier modes are discrete, and allow a 
clear separation between the non-oscillating part (the constant mode) and
the oscillating ones (the non-constant modes). But in the non-periodic
case, there is no such separation, and Fourier modes close to zero
create trouble. To control these low frequencies, we must 
  again inject some probabilistic information (namely, the ergodic
  theorem).  This difficulty
appears also further on in the study, to establish energy estimates. 

The rest of the paper is structured as follows.  The next section contains
a precise modeling of the domain, and the statements of the mathematical
results. The third section is devoted to the Dirichlet wall law. The fourth
section focuses on the boundary layer analysis. The final section is the
justification of Navier's wall law.

%%%%%%%%%%%%%%%%%%%%%%%%%%%%%%%%%%%%%%%%%%%%%%%%%%%%%%%%%%%%%%%%%%%%%%%%%%%
%%%%%%%%%%%%%%%%%%%%%%%%%%%%%%%%%%%%%%%%%%%%%%%%%%%%%%%%%%%%%%%%%%%%%%%%%%%  
\section{Statement of the results}

%%%%%%%%%%%%%%%%%%%%%%%%%%%%%%%%%%%%%%%%%%%%%%%%%%%%%%%%%%%%%%%%%%%%%%%%%%%
\subsection{Modeling of the rough domain} \label{modeling}

Let $\eps > 0$, and $(M, {\cal M}, \mu)$ a probability space. 
For all $m \in M$, we define
 a rough domain $\Omega^\eps(m)$ by 
$$ \Omega^\eps(m)  \: = \: \Omega \cup R^\eps_l(m) \cup R^\eps_u(m) , $$ 
where $\Omega = \R \times ]0,1[$ is the interior domain, and
    $R^\eps_{l.u}(m)$  is the lower, upper rough part. To obtain a
    realistic model for roughness,  
    we use  spatially homogeneous random fields:  following 
\cite{Papanicolaou:1982}, \cite{Jikov:1994} or
\cite{Bourgeat:1994},
 we recall that a homogeneous random field   is a measurable map 
$$\gamma : M \times \R^n \mapsto \R^m$$ 
satisfying: for all $h, z_1, \dots, z_k$ and all Borel
subsets $B_1, \dots, B_k$ of $\R^m$, 
\begin{multline*}
\mu\left( \left\{ m \in M, \: \gamma(m, z_1+h) \in B_1, \dots,  
\gamma(m, z_k+h) \in B_k \right\} \right) \\ 
= \mu\left( \left\{ m \in M, \:
\gamma(m, z_1) \in B_1, \dots, \gamma(m, z_k) \in B_k \right\} \right). 
\end{multline*}
We remind that for $n=m=1$, a homogeneous random field is often called a
stationary random process. We thus  define 
\begin{align*}
 R^\eps_l(m) & = \left\{ x = (x_1,x_2), \: x_1 \in \R^2, \: 0 > x_2 > - \eps
 \gamma_l(m, x_1/\eps) \right\}, \\
 R^\eps_u(m) & = \left\{ x = (x_1,x_2), \: x_1 \in \R^2, \: 0 < x_2 - 1  <  \eps
 \gamma_u(m, x_1/\eps) \right\}, 
\end{align*}
where  $(\gamma_l, \gamma_u) = (\gamma_l, \gamma_u)(m, y_1) \in ]0,1[^2$
 is a homogeneous random
field. Moreover, we assume that for all $m$, $\gamma_{l,u}(m, \cdot)$ is
 a K-Lipschitz function, with $K > 0$ independent of $m$, and that $m
 \mapsto \gamma_{l,u}(m, \cdot)$ is measurable with values in the set
 $C_b(\R; \, \R^2)$ of continuous bounded functions.

Following a classical construction of Doob (see \cite{Doob:1990}  or 
\cite{Bourgeat:1994} for all necessary details), we 
introduce  another probability space, which will be more  convenient to our
description. Let $P$ the set of K-Lipschitz functions $\omega   : \R
\mapsto ]0,1[^2$. Let ${\cal P}$ the $\sigma$-algebra generated by the sets 
$$ \left\{ \omega \in P, \: \omega(y_1) \in A, \: \forall y_1 \in B \right\}, $$
 where $B$ is a finite subset of $\Q$,
 and $A$ is a disk with rational center and radius. Note that ${\cal P}$ is
 simply the borelian $\sigma-$algebra of $P$, seen as a  subset of 
 $C_b(\R; \, \R^2)$. Finally, consider
 the set function $\pi : {\cal P} \mapsto \R$ given by 
 $$ \pi(Q) = \mu\left( \left\{ m \in M, \: (\gamma_u, \gamma_l)(m, \cdot)
 \in Q \right\}\right). $$
One can show ({\it cf} \cite{Doob:1990}) 
 that $\pi$ is a probability measure on $(P, {\cal P})$. 
Moreover, we can define a translation group 
$$\tau_{h}: P \mapsto P, \quad \tau_h(\omega)(y_1) = \omega(y_1+h), \quad \forall y_1,
h \in \R,$$
 that preserves $\pi$.
In this way, one can also  describe the boundaries 
 with the measurable  map 
$$ (\omega_l, \omega_u)(y_1, \omega) = \omega(y_1). $$
Indeed,  the  laws of the random variables 
$$ \left((\gamma_l, \gamma_u)(\cdot,z_1+h), \dots, (\gamma_l,
\gamma_u)(\cdot,z_k+h)\right) $$  
and
$$ \left((\omega_l, \omega_u)(\cdot,z_1+h),  \dots, (\omega_l,
\omega_u)(\cdot,z_k+h)\right) $$  
are the same (and independent of $h$, for any $z_1, \dots, z_k$). 
{\em The advantage of this last framework is that one can write 
$$\omega_{u,l}(y_1, \omega) = h_{u,l} \circ \tau_{y_1}(\omega), 
\mbox{ with } (h_u, h_l)(\omega) = \omega(0), $$ 
 which will be useful in the study of boundary layers. 
 Hence, we will rather use the  
 formulation in terms of $\omega, h_u, h_l$}, and consider, for
all $\omega \in P$,   
\begin{equation}
 \Omega^\eps(\omega),  \: = \: \Omega \cup R^\eps_l(\omega) \cup
 R^\eps_u(\omega) ,
\end{equation}
where $\Omega = \R \times ]0,1[$, 
 \begin{equation} 
\begin{aligned}
 R^\eps_l(\omega) & = \left\{ x = (x_1,x_2), \: 0 > x_2 > - \eps
 h_l \circ \tau_{x_1/\eps}(\omega) \right\}, \\
 R^\eps_u(\omega) & = \left\{ x = (x_1,x_2), \: 0 < x_2 - 1  <
 \eps  h_u \circ \tau_{x_1/\eps}(\omega)\right\}. 
\end{aligned}
\end{equation}
We also define $\Sigma_0 = \R \times \{ 0 \}$, $\Sigma_1 = \R \times
\{1\}$ the horizontal boundaries of  $\Omega$.
\begin{figure}
\begin{center}
\includegraphics[height = 6.5cm, width=7.5cm]{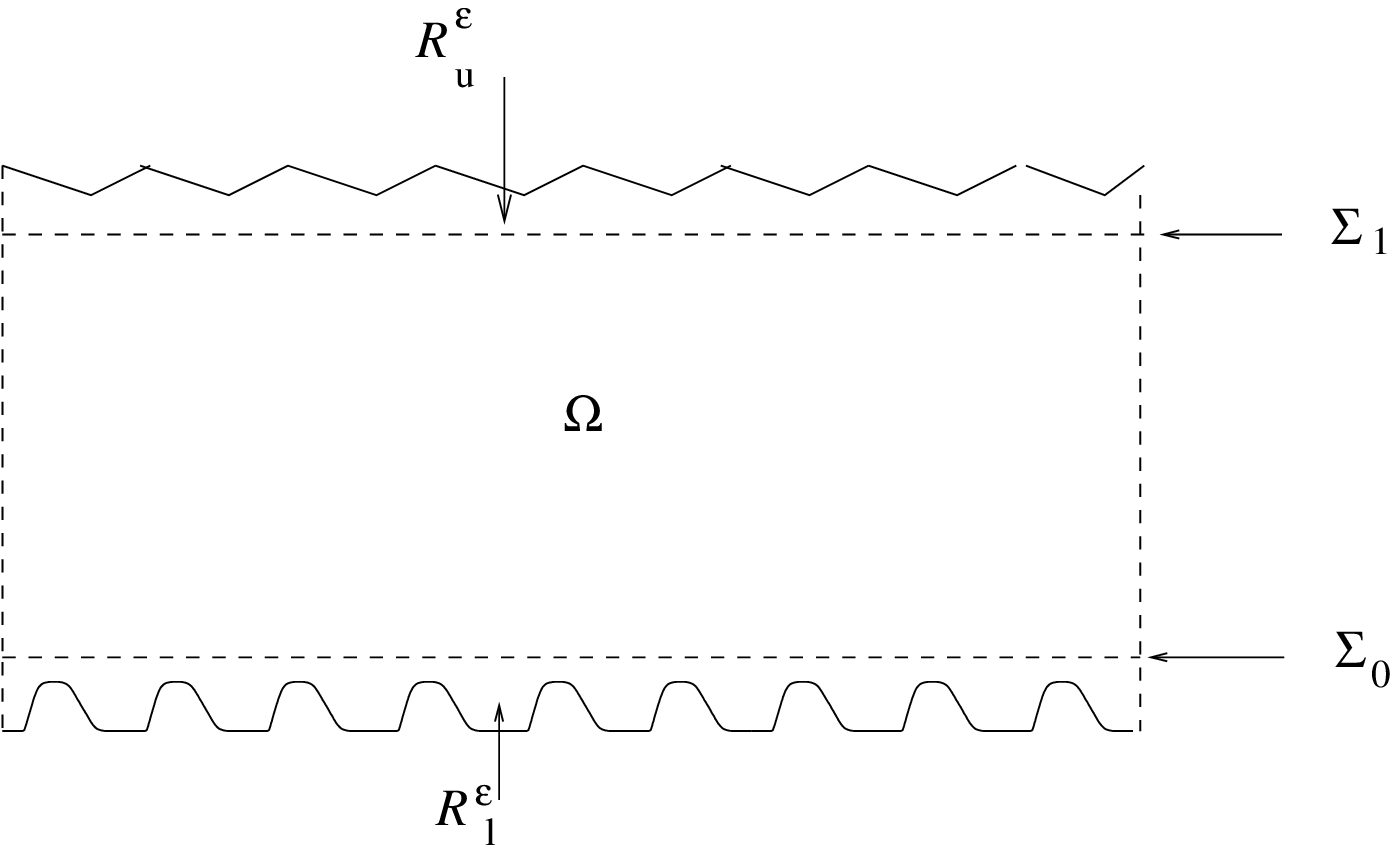}
\end{center}
\caption{The rough domain $\Omega^\eps$.}
\end{figure}

Note that our  modeling is derived from stochastic  
 homogeneization, where domains with small holes are described with such
  random fields (see again \cite{Papanicolaou:1982},  
\cite{Jikov:1994}, \cite{Chabi:1995}, \cite{Abddaimi:1996}, 
\cite{Beliaev:1996} among others). However, we emphasize that the classical
 tools of homogeneization 
(such as compensated compactness, two-scale convergence) do not apply to 
our boundary problem: broadly,  boundary layers are not seen in weak 
convergence processes, so that they do not allow to recover precise energy
estimates. To prove the theorems of the next sections will require a
precise construction of approximate solutions.   

Note also  that this random framework includes the periodic one. Take
$P = \T$ the unit torus, ${\cal P}$ its borelian $\sigma$-algebra, and
$\pi$ the Lebesgue measure. For any periodic function $F$ on $\T$,
any $y_1,\omega$ in $\T$ the formula 
$$ f(y_1,\omega) = F(y_1 + \omega) = F\left( \tau_{y_1}(\omega)\right) $$
allows to ``randomize'' the periodic structure.

In addition to the rough domain $\Omega^\eps$, we need to define  {\em
  boundary layer domains}, that will be useful in the study of Navier's
wall law. Namely, 
\begin{equation}
\begin{aligned} 
& {\cal R}_l(\omega) = \left\{ y = (y_1,y_2), \: y_2 > 
h_l \circ \tau_{y_1}(\omega)
\right\}, \\
& {\cal R}_u(\omega) = \left\{ y = (y_1,y_2), \: y_2 < h_u \circ
\tau_{y_1}(\omega)  \right\}.
\end{aligned}
\end{equation}
We will denote 
$$ {\cal R}_{l,u}(\omega) =  {\cal R}_{l,u}^+(\omega) \cup \Sigma_0 \cup
{\cal R}_{l,u}^-(\omega), $$ 
with 
$$ {\cal R}_{l,u}^{\pm}(\omega) = {\cal R}_{l,u}(\omega) \cap \{ \pm y_2 >
0 \}. $$

%%%%%%%%%%%%%%%%%%%%%%%%%%%%%%%%%%%%%%%%%%%%%%%%%%%%%%%%%%%%%%%%%%%%%%%%%%%
\subsection{Main results}

%%%%%%%%%%%%%%%%%%%%%%%%%%%%%%%%%%%%%%%%%%%%%%%%%%%%%%
\subsubsection{Dirichlet wall law} \label{Dwalllaw}

The first step in the analysis of system \eqref{NS}, \eqref{BC}, \eqref{flux}
is to prove the existence and uniqueness of solutions.
 For  given parameters $\eps, \omega$, this question has been adressed by 
Lady\v{z}enskaja and Solonnikov in article \cite{Ladyvzenskaja:1983}. We
also refer to \cite{Galdi:1994} for good overview on channel flow problems. 
The results of \cite{Ladyvzenskaja:1983} yield 
 existence and uniqueness of a solution $u^\eps(\omega, \cdot)$ in the space
$$ B_2(\Omega^\eps(\omega)) = 
\left\{ u \in H^1_{0,loc}(\Omega^\eps(\omega)), \quad  \sup_{R \ge 1}
\frac{1}{R} \int_{\Omega^\eps(\omega, R)} | \na u |^2  < +\infty
\right\},  
$$  
where 
$$ \Omega^\eps(\omega, R) = \Omega^\eps(\omega) \cap \left\{ 0 < |x_1| < R
  \right\}, \quad R >0, $$
{\em for small enough flux  $\phi < \phi_0$}.
 However, we can not apply  this 
result, as the  dependence of  $\phi_0$ with respect to $\eps$ and
$\omega$, as well as measurability properties with respect to $\omega$ are
not clear. Hence, we  show directly the following existence and uniqueness
result:  
\begin{theorem} \label{existenceNS}
There exists $\phi_0 > 0$ and  $\eps_0 > 0$  such that: 
for all $\phi < \phi_0$, for all $\eps < \eps_0$, for all $\omega \in P$,
system \eqref{NS}, \eqref{BC}, \eqref{flux} has a unique solution
$$u^\eps(\omega, \cdot) \in B_2\left(\Omega^\eps(\omega)\right). $$
Moreover, if we  denote  $\tilde{u}^\eps(\omega, \cdot)$ the zero-extension
  of $u^\eps(\omega, \cdot)$ outside $\Omega^\eps(\omega)$, then the
mapping
$$ P \mapsto H^1_{loc}(\R^2), \quad \omega \mapsto \tilde{u}^\eps(\omega,
\cdot) $$
is measurable. 
\end{theorem}
As will be clear from the proof, the solution $u^\eps$ is a perturbation of
the following flow:
\begin{equation} \label{poiseuille}
\left\{
\begin{aligned}
u^0(\omega, x) 
 & = \left(6\phi x_2 (1 - x_2), 0\right),  \: x \in \Omega\\
u^0(\omega, x)  & = 0, \: x \in \Omega^\eps-\Omega.
\end{aligned}
\right.
\end{equation}
Note that in restriction to $\Omega$, $u^0$ is simply the {\em Poiseuille
  flow}, that is the solution of \eqref{reduced}, \eqref{Dirichlet}. We
show the following estimates
\begin{theorem} \label{ThDirichlet}
There exists $C > 0$ such that, for all $\omega \in P$:
\begin{equation} \label{estimatesDirichlet}
\left\{
\begin{aligned}
&\sup_{R \ge 1} \frac{1}{R} \int_{\Omega^\eps(\omega,R)} \left|\na u^\eps(\omega,
\cdot) - \na u^0(\omega, \cdot)\right|^2 \: \le \:  C \eps, \\
& \sup_{R \ge 1} \frac{1}{R} \int_{\Sigma_0(R) \cup \Sigma_1(R)}
 \left| u^\eps(\omega, \cdot) \right|^2  \: \le 
\:  C \eps^2, \\
&\sup_{R \ge 1} \frac{1}{R} \int_{\Omega(R)} \left|  u^\eps(\omega, \cdot) 
- u^0(\omega, \cdot) \right|^2  \: \le  \:  C \eps^2,
\end{aligned}
\right.
\end{equation} 
where 
$$\Omega(R) =  \Omega \cap \left\{0 < |x_1| < R\right\}, \quad \Sigma_{0,1}(R) =
\Sigma_{0,1}  \cap \left\{ 0 < |x_1| <  R\right\}. $$
\end{theorem} 
Theorem \ref{estimatesDirichlet} expresses that the Dirichlet wall law
yields a $O(\eps)$ quadratic approximation of the real solution. The obtention
of the $L^2$-estimate relies on a duality argument.

%%%%%%%%%%%%%%%%%%%%%%%%%%%%%%%%%%%%%%%%%%%%%%%%%%%%%%
\subsubsection{Navier's wall law} \label{Nwalllaw}

The Dirichlet wall law, that is  the approximation of $u^\eps$ by $u^0$
does not account for the
behaviour of $u^\eps$ near the rough bondaries. To derive a more accurate
wall law, we carry a boundary layer analysis.  We show that, for
small $\eps$,  $u^\eps$ is close to 
 \begin{equation} \label{Ansatz} 
\begin{aligned}
u^\eps_{app}(\omega, x) & = u^0(\omega, x)  + \eps  U_l\left(\omega, 
\frac{x_1}{\eps}, \frac{x_2}{\eps} \right) + \eps 
 U_u\left(\omega, \frac{x_1}{\eps},
\frac{x_2 - 1}{\eps} \right)  \\
&\quad + \eps  u^1(\omega,x),
\end{aligned}  
\end{equation}
where $U_{u,l}$ is an appropriate boundary layer term, and $u^1$ an
appropriate Poiseuille type flow.  Precisely, $U_{l,u} = U_{l.u}(\omega, y)$ 
satisfies a Stokes system with jump conditions:
\begin{equation} \label{BL}
\left\{
\begin{aligned}
&-\nu \Delta U_{l,u}(\omega, \cdot) + \na P_{l,u}(\omega, \cdot) = 0, 
\quad  y \in {\cal R}_{l,u}^\pm(\omega) \\
&\div U_{l,u}(\omega, \cdot) = 0, \quad y \in {\cal R}_{l,u}(\omega),\\
&\left[ U_{l,u}(\omega,\cdot) \right]\vert_{\Sigma_0} = 0, \\
&\left[ \pa_2 U_{l,u}(\omega, \cdot) -  P_{l,u}(\omega,
  \cdot) e_2 \right]\vert_{\Sigma_0} = -6\phi, \\
& U_{l,u} = 0, \quad y \in \pa {\cal R}_{l,u}(\omega).
\end{aligned}
\right.
\end{equation}
In the case of periodic roughness, the solvability of system \eqref{BL} is
direct. Morever, using the Fourier transform in tangential variables, it can
be seen easily that {\em $U_{l,u}$ converges exponentially fast to a
  constant} as $y_2$ goes to infinity.
 The random setting requires much more work. We prove in
section \ref{sectionBL} the following 
\begin{theorem} \label{existenceBL}
For all $\phi > 0$, there exits a unique variational solution $U_{l,u}$
 of \eqref{BL}, in the sense of paragraph \ref{defvariational}.
 It is almost surely a classical solution, and satisfies 
\begin{equation}  \label{boundsBL}
\begin{aligned}
& \sup_{R \ge 1} \frac{1}{R} \, \E \left(\int_{{\cal R}_{l,u}(\cdot,R)} \left| \na
 U_{l,u}(\cdot, y) \right|^2 dy     \right) < +\infty, \\
& \sup_{R \ge 1} \frac{1}{R} \int_{{\cal R}_{l,u}(\omega,R)}  \left| \na
 U_{l,u}(\omega, y) \right|^2 dy < \infty \quad \mbox{ almost surely },
\end{aligned}
\end{equation}  
where as usual 
$$\displaystyle 
{\cal R}_{l,u}(\omega, R) =  {\cal R}_{l,u}(\omega) \cap \{ 0 < |y_1| < R
\}. $$
  Moreover, there exists a measurable map: $\displaystyle 
U^{\infty}_{l,u}: P \mapsto \R^2$, with $U^\infty_{l,u,2} = 0$, such that 
\begin{equation} \label{cvgce_BL}
\begin{aligned}
& \sup_{R \ge 1}   \, \E \left(  \frac{1}{R} \int_{|y_1| < R}
  |U_{l,u}(\omega, y_1, y_2) dy_1 -
U^{\infty}_{l,u}(\omega) |^2 \,dy_1  \right) \xrightarrow[y_2 \rightarrow
  \infty]{}  0, \\
&   \left| U_{l,u}(\omega, y_1, y_2) - U^{\infty}_{l,u}(\omega) \right|
 \xrightarrow[y_2 \rightarrow
  \infty]{}  0, \: \mbox{ locally uniformly in } y_1.
\end{aligned}
\end{equation} 
\end{theorem}
The convergence of $U_{l,u}$ towards a constant is proved using the
ergodic theorem. Contrary to the periodic case, 
we are not able to precise the speed of convergence (see section
\ref{sectionBL}  for all details).

We then establish the  estimates on $u^\eps - u^\eps_{app}$:
\begin{theorem} \label{thNavier}
Let $U_{l,u}$ as in theorem \ref{existenceBL}, $u^1$ given in \eqref{u1}.
  The approximation $u^\eps_{app}$ of \eqref{Ansatz} satisfies, as
$\eps \rightarrow 0$:
\begin{equation} \label{estimatesNavier}
\left\{
\begin{aligned}
&\sup_{R \ge 1} \frac{1}{R} \, \E \left( \int_{\Omega^\eps(\cdot,R)}
 \left|\na u^\eps(\omega,
\cdot) - \na u^\eps_{app}\right|^2 \right) \: = \: o(\eps^2), \\
%& \sup_{R \ge 1} \frac{1}{R}  \, \E \left( \int_{\Sigma_0(R) \cup \Sigma_1(R)}
% \left| u^\eps(\omega, \cdot) - u^\eps_{app}(\omega, \cdot)\right|^2
% \right) \: = 
%\:  o(\eps^3), \\
&\sup_{R \ge 1} \frac{1}{R} \, \E \left( 
\int_{\Omega(R)} \left|  u^\eps(\omega, \cdot) 
- u^\eps_{app}(\omega,\cdot) \right|^2  \right)  \: =  \:  o(\eps^2),
\end{aligned}
\right.
\end{equation} 
\end{theorem}
In the periodic setting, for which the boundary layer profiles
 $\displaystyle 
U_{l,u} - U^\infty_{l,u}$ decay exponentially with $y_2$, the bound
$o(\eps^2)$  turns to  $\displaystyle 
O(\eps^3)$ for the $H^1$ estimate,
 and $\displaystyle O(\eps^4)$  for the $L^2$ estimate.
 But in the general random setting, we are not able to improve  
the bound $o(\eps^2)$: 
it is  related to the speed of convergence of
 $U_{l,u}$ as $y_2 \rightarrow +\infty$.
 Theorem \ref{thNavier}  extends  the results of J\"ager and Mikelic,
 \cite[theorem 1, p113]{Jager:2001} for periodic roughness.
 As they consider domains  with lateral boundaries,
 they need  to add a lateral boundary layer
 term, which yields a less  precise estimate: the $O(\eps^3)$ and
 $O(\eps^4)$ bounds are replaced by $O(\eps^2)$ and $O(\eps^3)$ bounds
 respectively. The unbounded setting allows to avoid this loss of
 accuracy.

We may now justify the relevance of Navier's friction law. Let
$v^\eps(\omega, \cdot)$ the solution of 
\begin{equation} \label{Naviereq}
     \left\{
\begin{aligned}
& v \cdot \na v + \na q - \nu  \Delta v = 0, \: x \in \Omega, \\
& \div v = 0, \: x \in \Omega, \\
& \int_{\sigma(x_1)} \!\!\!\! v_1 \:= \: \phi,  \quad v^\eps_2(\omega,
  \cdot)\vert_{\Sigma_{0,1}} = 0,
\end{aligned}
\right.
\end{equation} 
 with Navier's condition 
\begin{equation} \label{Navierlaw}
\begin{aligned}
& v^\eps_1(\omega, \cdot) \: = \:  + \eps \, \alpha_l(\omega) \,    \, 
\frac{\pa v^\eps_1(\omega, \cdot)}{\pa x_2},  \quad x_2 = 0, \quad
\alpha_l(\omega) 
= \frac{U_{l,1}^\infty(\omega)}{6 \phi}, \\
& v^\eps_1(\omega, \cdot) \: = \:  -\eps \alpha_u(\omega)  \, \frac{\pa
    v^\eps_1(\omega, \cdot)}{\pa x_2},  \quad x_2 = 1, \quad \alpha_u(\omega)
= \frac{U_{u,1}^\infty}{6 \phi},
\end{aligned}
\end{equation}
 with $U^\infty_l, U^\infty_u$ as in theorem \ref{existenceBL}. Note that
 by linearity of \eqref{BL}, {\em constants $\alpha_{l,u}(\omega)$ depend only
 on the roughness, not on the flux $\phi$.}  We
state:   
\begin{theorem} \label{justifNavier}
$$ \sup_{R \ge 1} \frac{1}{R} \, \E \left( 
\int_{\Omega(R)} \left|  u^\eps(\omega, \cdot) - v^\eps(\omega, \cdot)
 \right|^2  \right)  \: = \:  o(\eps^2), \quad \eps \rightarrow 0. $$
\end{theorem}
Thus, the Navier's friction law leads to a (slightly) better approximation
than Dirichlet's law. Again, the proof of theorem \ref{justifNavier}
involves the ergodic theorem, which prevents  quantitative bounds
(the $o(\eps)$ can not {\em a priori} be precised). 
We end this presentation of the results with two remarks: 
\begin{enumerate}
\item All the estimates of section \ref{Nwalllaw} involve  the expectation
  of the spatial quadratic norms. Due to a lack of deterministic control
  of the boundary layers, we are not able to obtain almost sure
  estimates in $B_2(\Omega^\eps)$,  like in theorem
  \ref{ThDirichlet}. However, {\em one can deduce informations 
with high probability}. For instance, theorem
  \ref{justifNavier} and Tchebitchev inequality show that: for all $R$, 
$\delta > 0$, 
$$ P\left( \omega, \: 
\int_{\Omega(R)} \left|  u^\eps(\omega, \cdot) - v^\eps(\omega, \cdot) 
\right|^2 > \delta \eps^2  \right) \: \xrightarrow[\eps \rightarrow 0]{}
0. $$
\item Under an assumption of ergodicity on the group $(\tau_h)$, that is 
$$ \forall  A \in {\cal P}, \:  \left( \tau_h(A) = A, \forall h \right)
  \: \Rightarrow \: \left( \pi(A) = 0
  \mbox{ or } \pi(A) = 1 \right),$$
 {\em the constants $\alpha_{l,u}$ that appear in \eqref{Navierlaw}
 are independent of $\omega$}, as a consequence of the ergodic theorem 
(see section \ref{sectionBL}).
 In such case, one does not need to know the shape of
 the boundary to determine the appropriate coefficient in Navier's wall law. 
 It can be deduced almost surely from a numerical computation  involving
 another boundary. 
\end{enumerate}

%%%%%%%%%%%%%%%%%%%%%%%%%%%%%%%%%%%%%%%%%%%%%%%%%%%%%%%%%%%%%%%%%%%%%%%%%%%%
%%%%%%%%%%%%%%%%%%%%%%%%%%%%%%%%%%%%%%%%%%%%%%%%%%%%%%%%%%%%%%%%%%%%%%%%%%%%
\section{Justification of Dirichlet wall law} \label{sectionDirichlet}

%%%%%%%%%%%%%%%%%%%%%%%%%%%%%%%%%%%%%%%%%%%%%%%%%%%%%%%%%%%%%%%%%%%%%%%%%%%%
\subsection{Well-posedness}

This section is devoted to the proof of theorem \ref{existenceNS}. The
solution $(u^\eps, p^\eps)$ is searched as a perturbation of the Poiseuille type
flow $(u^0, p^0)$, where $u^0$ is defined by  \eqref{poiseuille}, and 
$$ p^0(\omega,x) = -12\nu \phi x_1, \quad \forall x \in
\Omega^\eps(\omega). $$ 
We denote $u^\eps = u^0 + w, \:$ $p^\eps = p^0 + q$, and consider the
following system (dependence on $\omega$ is omitted to lighten
notations): 
\begin{equation} \label{NS2}
\left\{
\begin{aligned} 
& w \cdot \na w + u^0 \cdot \na w + w \cdot \na u^0 + \na q - \nu  \Delta w
  = f^\eps,  \quad  x \in \Omega^\eps\setminus\left(\Sigma_0 \cup \Sigma_1\right), \\
& \div w = 0, \quad  x \in \Omega^\eps, \\
&  [w]\vert_{\Sigma_0 \cup \Sigma_1}
  =0, \quad \left[ \pa_2 w - q e_2 \right]\vert_{\Sigma_0
    \cup \Sigma_1} = -
 6\phi\\
&  \int_{\sigma(x_1)} \!\!\!\! w_1 \:= \: 0, \quad  w\vert_{\pa \Omega^\eps} = 0,
\end{aligned}
\right.
\end{equation} 
where 
\begin{align*}
 f^\eps(\omega, x) & = 0, \quad x \mbox{ in } \Omega, \\
 f^\eps(\omega, x) & = \left( -12\nu \phi,
0 \right), \quad x \mbox{ in } \Omega^\eps(\omega) \setminus \Omega. 
\end{align*}
The proof of theorem \ref{existenceNS} is divided in three steps:
\begin{enumerate}
\item We consider  the linear problem
\begin{equation}  \label{linearNS}
\left\{
\begin{aligned} 
& (v + u^0)
 \cdot \na w + w \cdot \na u^0 + \na q - \nu  \Delta w
  = f + \div G, \\
& x \in \Omega^\eps\setminus\left(\Sigma_0 \cup \Sigma_1\right), \quad
 \div w = 0, \quad  x \in \Omega^\eps, \\
&  [w]\vert_{\Sigma_0 \cup \Sigma_1}
  =0, \quad \left[ \pa_2 w - q e_2 \right]\vert_{\Sigma_0
    \cup \Sigma_1} = \tilde{\phi}\\
&  \int_{\sigma(x_1)} \!\!\!\! w_1 \:= \: 0, \quad  w\vert_{\pa \Omega^\eps} = 0,
\end{aligned}
\right.
\end{equation} 
  where for all $\omega$, 
\begin{align*} 
& G(\omega, \cdot) \in L^2_{uloc}\left(\Omega^\eps(\omega)\right) =
 \Bigl\{ G, \: \sup_{\tau > 0} \int_{\substack{x \in
     \Omega^{\eps}(\omega)\\
     \tau < |x_1| < \tau + 1}} \!\!\!\!  |G |^2 < +\infty \Bigr\},\\
& f(\omega, \cdot) \in L^2_{uloc}\left(\Omega^\eps(\omega)\right), \: f = 0
 \mbox{ in } \Omega, \\
&v(\omega, \cdot) \in H^1_{uloc}\left(\Omega^\eps(\omega)\right) = 
\Bigl\{ v, \quad (v, \na v) \in
L^2_{uloc}\left(\Omega^{\eps}(\omega)\right)\Bigr\}, \\
& \tilde{\phi}(\omega, \cdot) \in L^2_{uloc}(\R). 
\end{align*}   
We show the existence and uniqueness of a solution  
$$ w(\omega, \cdot) \in H^1_{uloc}\left(\Omega^\eps(\omega)\right), \quad
\forall \omega. $$
\item  Thanks to this linear analysis, we show that, 
for $\eps<\eps_0$, $\phi<\phi_0$  small enough, there exists a unique 
  solution $w$ of \eqref{NS}, thus of a solution $u^\eps$ of \eqref{NS}, as
  in theorem \ref{existenceNS}.
\item We show measurability properties of $u^\eps$. 
\end{enumerate}

%%%%%%%%%%%%%%%%%%%%%%%%%%%%%%%%%%%%%%%%%%%%%%%
\subsubsection{Linear problem}

To build a solution of \eqref{linearNS}, we consider the following
approximate problem:
\begin{equation} \label{linearNS2}
\left\{
\begin{aligned} 
& (v + u^0)
 \cdot \na w_n + w_n \cdot \na u^0 + \na q_n - \nu  \Delta w_n
  = f + \div G, \\
& x \in \Omega^\eps_n\setminus\left(\Sigma_0(n) \cup \Sigma_1(n)\right), \quad
 \div w_n = 0, \quad  x \in \Omega^\eps_n, \\
&  [w_n]\vert_{\Sigma_0(n) \cup \Sigma_1(n)}
  =0, \quad \left[ \pa_2 w_n - q e_2 \right]\vert_{\Sigma_0(n)
    \cup \Sigma_1(n)} = \tilde{\phi}\\
& w_n\vert_{\pa \Omega^\eps_n} = 0,
\end{aligned}
\right.
\end{equation} 
where $\Omega^\eps_n$ is a short-hand for the bounded domain 
$\Omega^\eps(\omega,n)$. A formal energy estimate yields
\begin{align*}
  \nu  \int_{\Omega^\eps_n} \left| \na w_n \right|^2   \: = \:  & 
 \int_{\Omega^\eps_n\setminus\Omega(n)} \!\!\!\!\!\!  \!\!\!
 f w_n  \, - \, \int_{\Omega^\eps_n} G \na w_n \,
 - \,   \int_{\Omega^\eps_n} (w_n \cdot \na u^0) w_n \\
&  \: + \: \tilde{\phi}  \int_{\Sigma_0(n)} \!\!\!\!\!\! w_n - 
\tilde{\phi} \int_{\Sigma_1(n)} \!\!\!\!\!\!  w_n  \\
 \: \le \:  &  C \biggl(  \sqrt{n}  \,  \| f \|_{L^2_{uloc}}  \: \| w_n
\|_{L^2(\Omega^\eps_n\setminus\Omega(n))} \\
& \: + \:  \, \sqrt{n} \,  \| G
\|_{L^2_{uloc}} \, \| \na w_n \|_{L^2(\Omega^\eps_n)} 
\\
&  \: + \:   \, \phi  \,  \| w_n
\|^2_{L^2(\Omega^\eps_n)} \: + \:    \, \| \tilde{\phi} \|_{L^2_{uloc}}
\sqrt{n} \:
  \| w_n \|_{L^2(\Sigma_0(n) \cup \Sigma_1(n))} \biggr).
\end{align*}
We remind the Poincar\'e inequality
\begin{equation} \label{Poincare}
\| \varphi \|_{L^2({\cal O}_a)} \: \le \: C \, a \, \| \na \varphi \|_{L^2({\cal
    O}_a)}
\end{equation}
valid in a domain 
\begin{equation} \label{Oa}
 {\cal O}_a = \{ x,  \: \alpha < x_1 < \beta,  \: \gamma_l(x_1) < x_2  <
\gamma_u(x_1) \}, \quad \sup | \gamma_l - \gamma_u | < a, 
\end{equation}
 for all $\varphi$ in $H^1({\cal O}_a)$, with $\varphi=0$ on one of the 
 boundaries $\gamma_{l,u}$. It is easily deduced from \eqref{Poincare}
 that: for all 
$\varphi \in H^1(\Omega^\eps_n\setminus\Omega(n))$, $\varphi = 0$ on $\pa
 \Omega^\eps_n\setminus\Omega(n)$,   
\begin{equation} \label{Poincare2}
\begin{aligned}
& \| \varphi \|_{L^2(\Omega^\eps_n - \Omega(n))} \: \le \:  C \,  \eps \, 
 \| \na \varphi
 \|_{L^2(\Omega^\eps_n - \Omega(n))}, \\
& \| \varphi \|_{L^2(\Sigma_0(n) \cup \Sigma_1(n))} \: \le \:  C  \,
 \eps^{1/2} \,
 \| \na \varphi  \|_{L^2(\Omega^\eps_n\setminus \Omega(n))}.
\end{aligned}
\end{equation}
We infer  from \eqref{Poincare}, \eqref{Poincare2}  that
\begin{align*}
 \int_{\Omega^\eps_n} \left| \na w_n \right|^2 \:\le \: & 
 C \, \Bigl( \sqrt{n}
\left(\eps \| f \|_{L^2_{uloc}} +  \| G
\|_{L^2_{uloc}} + \sqrt{\eps}  \| \tilde{\phi} \|_{L^2_{uloc}} \right) 
    \| \na w_n \|_{L^2(\Omega^\eps_n)} \\
& \:  + \:   \phi \, \| \na w_n \|_{L^2(\Omega^\eps_n)}^2  \Bigr)
\end{align*}
Hence, 
for $\phi < \phi_0$, $\eps < \eps_0$ with $\phi_0, \: \eps_0$ small
enough, we get (uniformly in the variable $\omega \in P$) 
\begin{equation} \label{estimwn}
\frac{1}{\sqrt{n}} \, \| \na w_n \|_{L^2(\Omega^\eps_n)} \: \le \:  C \,
\left(  \eps \| f \|_{L^2_{uloc}} \, + \,  \| G
\|_{L^2_{uloc}} + \sqrt{\eps} \| \tilde{\phi} \|_{L^2_{uloc}} \right).
\end{equation}
By Lax-Milgram lemma, this estimate yields existence and uniqueness of a
solution $w_n$ of \eqref{linearNS2}. 

We now wish to let $n$ go to infinity. However, estimate \eqref{estimwn}
degenerates in this limit. To obtain compactness, we will follow ideas of    
 \cite{Ladyvzenskaja:1983} and localize estimate \eqref{estimwn} in a band
 of width $\eta$. To do so,
 we multiply the first equation of $\eqref{linearNS2}$ by  $w_n$, and
 integrate over 
$$ \Omega^\eps_n(R) = \Omega^\eps_n \cap \{ 0 < |x_1| < R \}. $$
Proceeding as above, we get, for $\phi < \phi_0$ small enough,
\begin{align*} 
  \int_{\Omega^\eps_n(R)} \nu | \na w_n |^2 \:  \le \; & 
 C \left( \eps^2 \| f \|_{L^2_{uloc}}^2 \, + \,  \| G\|_{L^2_{uloc}}^2 + 
\eps \| \tilde{\phi} \|^2_{L^2_{uloc}} \right)  (R + 1)    \\
& \: + \left| \int_{|x_1| = R} \nu \pa_1 w_n \cdot w_n \right| \, + \, 
\left| \int_{|x_1| = R}  q_n w_{n,1} \right| \\
&  \: + \: \left| \int_{|x_1| = R} v_1 | u |^2  \right| + \left| 
\int_{|x_1| = R} G e_1 \cdot u  \right|.
\end{align*}
We then integrate with respect to $R$, fro $\eta$ to $\eta+1$. We deduce 
\begin{align*}
F(\eta) \: \le \:  &  C  \left(  \eps^2 \| f \|_{L^2_{uloc}}^2 \, + \, 
 \| G\|_{L^2_{uloc}}^2 + \eps
\| \tilde{\phi} \|_{L^2_{uloc}}^2 \right) \left( \eta + 1 \right)   \\ 
& \: + \left| \int_{\Omega^\eps_n(\eta, \eta+1)}  \pa_1 w_n \cdot w_n
 \right| + \: 
\left| \int_{\Omega^\eps_n(\eta, \eta+1)}  q_n w_{n,1} \right| \\
& \: + \:
 \left|\int_{\Omega^\eps_n(\eta, \eta+1)}  v_1 | w_n |^2  \right| + \left| 
\int_{\Omega^\eps_n(\eta, \eta+1)}  G e_1 \cdot  w_n  \right|.
\end{align*}
where 
$$ F(\eta) = \int_\eta^{\eta+1} \nu \int_{\Omega^\eps_n(R)} | \na w_n |^2 \,
dR, \quad \Omega^\eps_n(\eta, \eta+1) = \Omega^\eps_n \cap \{ \eta < | x_1 |
< \eta + 1 \}. $$
By Poincar\'e inequality \eqref{Poincare},  and Sobolev inequality 
\begin{equation} \label{Sobolev}
\| \phi \|_{L^4(O_a)} \le C \| \phi \|_{H^1(O_a)}    
\end{equation}
one has easily that 
\begin{align*}
& \left| \int_{\Omega^\eps_n(\eta, \eta+1)}  \pa_1 w_n \cdot w_n
 \right| \: \le \: C \, \int_{\Omega^\eps_n(\eta, \eta+1)} |\na w_n |^2, \\
& \left|\int_{\Omega^\eps_n(\eta, \eta+1)}  v_1 | w_n |^2  \right| \: \le \: 
C \,  \| v \|_{H^1_{uloc}(\Omega^\eps_n)} \, 
\int_{\Omega^\eps_n(\eta, \eta+1)} |\na w_n |^2, \\
&  \left| \int_{\Omega^\eps_n(\eta, \eta+1)}  G e_1 \cdot w_n  \right| \: \le \:
C \, \| G \|_{L^2_{uloc}(\Omega^\eps_n)} \,  \left(
\int_{\Omega^\eps_n(\eta, \eta+1)} |\na w_n |^2 \right)^{1/2}. 
\end{align*}
The treatment of the integral involving the pressure is exactly the same as
in \cite{Ladyvzenskaja:1983}. Boundary conditions and incompressibility of
$w_n$ imply that $\int_{|x_1|= R} w_{n,1} = 0$ for all $R$, which yields 
$\int_{\Omega^\eps_n(\eta, \eta+1)} w_{n,1} = 0$. It is then well-known (see
\cite{Galdi:1994} for detailed description) that: there exists 
$$\varphi_\eta \in H^1_0(\Omega^\eps_n(\eta, \eta+1)),  \: \mbox{ with } 
\div \varphi_\eta = w_1, $$ 
and
$$ \| \varphi_\eta \|_{H^1_0(\Omega^\eps_n(\eta, \eta+1))} \: \le \:  C \| w_1
\|_{L^2(\Omega^\eps_n(\eta, \eta+1))}, $$
where $C$ is a positive constant independent on $\eps$, $n$ and $\eta$. We
can write 
$$  \int_{\Omega^\eps_n(\eta, \eta+1)}  q_n w_{n,1} =
\int_{\Omega^\eps_n(\eta, \eta+1)} q_n  \div \varphi_\eta = - 
  \int_{\Omega^\eps_n(\eta, \eta+1)} \na q_n \cdot \varphi_\eta. $$
Using the expression of $\na q_n$ in \eqref{linearNS2}, and after several
integrations by parts, 
\begin{align*}
 \left| \int_{\Omega^\eps_n(\eta, \eta+1)}  q_n w_{n,1} \right|
\: \le \:&  \left| \int_{\Omega^\eps_n(\eta, \eta+1)} \!\!\!\!\!\! \nu
\na w_n \cdot \na \varphi_\eta \right| \: + \: 
 \left| \int_{\Omega^\eps_n(\eta, \eta+1)} \!\!\!\!\!\!(v + u^0) \cdot \na w_n
 \varphi_\eta \right| \\
& +  \left| \int_{\Omega^\eps_n(\eta, \eta+1)}  \na \varphi_\eta w_n \cdot u^0
 \right| + \left| \int_{\Omega^\eps_n(\eta, \eta+1)} \!\!\!\!\!\! (f + \div G) w_n
 \right|  
\end{align*}
This leads to the inequality 
\begin{align*}
 \left| \int_{\Omega^\eps_n(\eta, \eta+1)}  q_n w_{n,1} \right| \: \le \: &  
C \left( 1 + \| v \|_{H^1_{uloc}} \right) 
 \int_{\Omega^\eps_n(\eta, \eta+1)} |\na w_n \|^2  \\
& \:+ C \left(  \eps^2 \| f \|_{L^2_{uloc}}^2 \, + \,  \| G
\|^2_{L^2_{uloc}}  \right). 
\end{align*}
Together with previous bounds, this yields 
\begin{equation}
\begin{aligned}
 F(\eta) \: \le  \: & C \left( 1 + \| v \|_{H^1_{uloc}} \right)
\int_{\Omega^\eps_n(\eta, \eta+1)} |\na w_n \|^2 \\
& \:  + \:  C
 \left(  \eps^2 \| f \|_{L^2_{uloc}}^2 \, + \,  \| G
\|_{L^2_{uloc}}^2 + \eps\| \tilde{\phi} \|_{L^2_{uloc}}^2 \right)  (1 + \eta). 
\end{aligned}
\end{equation}
that is
\begin{equation} \label{Gronwall}
 F(\eta) \le {\cal C}_1 F'(\eta) + \:  {\cal C}_2 
   (1 + \eta), 
\end{equation}
where 
$${\cal C}_1 = C_1  \left( 1 + \| v \|_{H^1_{uloc}} \right), \quad  {\cal
  C}_2  = C_2 \left(  \eps^2 \| f \|_{L^2_{uloc}}^2 \, + \,  \| G
\|_{L^2_{uloc}}^2 + \eps \| \tilde{\phi} \|_{L^2_{uloc}}^2 \right). $$
This last equation is a reverse Gronwall type inequality. By estimate
\eqref{estimwn}, up to take a larger $C_2$,
  we can suppose that 
$$ F(n) \: \le \: {\cal C}_2 \, ({\cal C}_1 + n + 1). $$ 
From \eqref{Gronwall}, we deduce easily that for all $\eta$,
\begin{align*}
& e^{-\frac{1}{{\cal C}_1} n} F(n) - e^{-\frac{1}{{\cal C}_1} \eta} F(\eta)  
 \ge \: - {\cal C}_2  
 \int_\eta^n \frac{(t + 1)}{{\cal C}_1} \, e^{-\frac{1}{{\cal C}_1} t} dt \\
& \ge \:  {\cal C}_2   ({\cal C}_1 + n + 1) e^{-\frac{1}{{\cal C}_1} n}   
-  {\cal C}_2 ({\cal C}_1 + \eta + 1) e^{-\frac{1}{{\cal C}_1} \eta}.
\end{align*} 
It follows that  for all $\eta$,
$$ F(\eta) \le {\cal C} \left( \eps^2 \| f \|_{L^2_{uloc}}^2 \, + \,  \| G
\|_{L^2_{uloc}}^2 + 
\eps \| \tilde{\phi} \|_{L^2_{uloc}}^2 \right) (\eta+1) $$
and finally 
\begin{equation} \label{B2estimatewn}
\sup_{\eta \ge 1} 
\frac{1}{\eta} \int_{\Omega^\eps_n(\eta)} |\na w_n |^2  \le {\cal C} 
\left( \eps^2 \| f \|_{L^2_{uloc}}^2 \, + \,  \| G
\|_{L^2_{uloc}}^2 + \eps \| \tilde{\phi} \|_{L^2_{uloc}}^2 \right),
\end{equation} 
where ${\cal C}$ can be taken affine in  $\| v \|_{H^1_{uloc}}$.

On the basis of such estimate, one can, for all $\omega$, extract a
subsequence $w_{\phi_\omega(n)}$ that converges weakly in 
$H^1_{loc}\left( \Omega^\eps(\omega)\right)$.  The
 limit $w(\omega, \cdot)$ is a solution of \eqref{linearNS}, and satisfies 
the estimate 
\begin{equation} \label{B2estimate}
\sup_{\eta \ge 1} 
\frac{1}{\eta} \int_{\Omega^\eps_n(\eta)} |\na w |^2  \le {\cal C} 
\left( \eps^2 \| f \|_{L^2_{uloc}}^2 \, + \,  \| G
\|_{L^2_{uloc}}^2 + \eps \| \tilde{\phi} \|_{L^2_{uloc}}^2 \right).
\end{equation} 

It remains to show that $w$ is in $H^1_{uloc}$. 
The argument is almost the 
same as in \cite[p745]{Ladyvzenskaja:1983}. 
Let $\tau > 3/2$. By energy estimates carried on the translated domain
$$\Omega^{\eps}(\tau, R) = \Omega^{\eps}(R) + (\tau,0),$$
we obtain similarly to \eqref{Gronwall}
\begin{equation} \label{Gronwall2}
 F_\tau(\eta) \le {\cal C}_1 F_\tau'(\eta) + \:  {\cal C}_2 
   (1 + \eta), \quad F_\tau(\eta) = \int_\eta^{\eta+1}
 \int_{\Omega^\eps(\tau, R)} | \na w |^2 \, dR. 
\end{equation}
By \eqref{B2estimate}, we have 
 $$F_\tau(\tau) \le {\cal C}  \left( \eps^2 \| f \|_{L^2_{uloc}}^2 \, + \,  \| G
\|_{L^2_{uloc}}^2 + \eps 
\| \tilde{\phi} \|_{L^2_{uloc}}^2 \right)  \, (\tau + 1),$$
 and reasoning as above, we get: for all $ 0 < \eta < \tau$,
$$ F_\tau(\eta) <  {\cal C}  \left( \eps^2 \| f \|_{L^2_{uloc}}^2 \, + \,  \| G
\|_{L^2_{uloc}}^2 + \eps \| \tilde{\phi} \|_{L^2_{uloc}}^2 \right)  \,
(\eta + 1).$$
Applying this inequality with $\eta = 3/2$ yields 
\begin{equation} \label{h1uloc}
\| w \|_{H^1_{uloc}}^2 \: \le \: {\cal C} 
\left( \eps^2 \| f \|_{L^2_{uloc}}^2 \, + \,  \| G
\|_{L^2_{uloc}}^2 + \eps \| \tilde{\phi} \|_{L^2_{uloc}}^2 \right)
\end{equation}
For uniqueness, one must show that the solution $w$ of \eqref{linearNS}
with $f=0$,  $G=0$, $\tilde{\phi}=0$ is identically zero. 
Inequality \eqref{Gronwall} turns into 
\begin{equation} \label{Gronwallunique}
 F(\eta) \: \le \:  {\cal C}_1 F'(\eta), \quad \forall \eta > 0, 
\end{equation}
from which it follows easily that 
$$ \limsup_{\eta \rightarrow +\infty} F(\eta) \, e^{-{\cal C}_1 \eta} \ge F(0).$$
Together with \eqref{B2estimate}, we deduce $w =0$. Note that uniqueness does
not only hold in $H^1_{uloc}(\Omega^\eps(\omega))$ but in the wider space
$B^2(\Omega^\eps(\omega))$, for any $\omega \in P$.

%%%%%%%%%%%%%%%%%%%%%%%%%%%%%%%%%%%%%%%%%%%%%% 
\subsubsection{Nonlinear problem}

This section is devoted to the well-posedness of \eqref{NS2}, for any
fixed $\omega \in P$. Therefore, we  consider linear equations \eqref{linearNS}
with special choice
$$ f = f^\eps, \: G = 0, \: \mbox{ and  } \tilde{\phi} =
-6\phi.$$
Note that $\| f^\eps \|_{L^2_{uloc}} \le C
\sqrt{\eps}$. We define the application 
$$ \Lambda :
H^1_{uloc}(\Omega^\eps(\omega)) \mapsto H^1_{uloc}(\Omega^\eps(\omega)),
\quad v(\omega, \cdot) \mapsto w(\omega, \cdot),$$ 
with $w(\omega, \cdot)$  the solution of \eqref{linearNS}. Then, for $\eps
< \eps_0$ small enough,  { $\Lambda$
\em is a contraction in restriction to the unit ball of 
$\displaystyle H^1_{uloc}(\Omega^\eps(\omega))$, }.

Indeed, let $v^1$, $v^2$ in the unit ball of 
$H^1_{uloc}(\Omega^\eps(\omega))$, and 
 denote $w^1 = \Lambda(v^1)$, $w^2 = \Lambda(v^2)$. By estimate
 \eqref{h1uloc}, they are bounded by 
$$ \| w^{1,2} \|_{H^1_{uloc}} \: \le \:  C \, \left( \eps^2 \| f^\eps
 \|_{L^2_{uloc}}^2  + \eps \| \phi \|_{L^2_{uloc}}^2 \right) \: \le \:  C
 \eps $$ 
for $\eps < \eps_0$ small enough. The difference $w = w^1- w^2$ satisfies 
\eqref{linearNS}, with 
$$ v = v^1, \: f = 0, \: G = (v^2 - v^1) \otimes w^2, \: \tilde{\phi} = 0. $$
 Again, using \eqref{h1uloc} leads to 
\begin{equation}
\begin{aligned}
\| w \|_{H^1_{uloc}} \: \le \: &  C \| (v^2 - v^1) \otimes w^2
 \|_{L^2_{uloc}} \: \le  C \: \| v^2 - v^1 \|_{H^1_{uloc}} \, \| w^2
 \|_{H^1_{uloc}} \\
 \le \:&  C \, \sqrt{\eps} \,  \| v^2 - v^1 \|_{H^1_{uloc}}.
\end{aligned}
\end{equation}
Hence, for $\eps < \eps_0$ small enough, $\Lambda$ is a  contraction in
restriction to the unit ball. 
By the Banach fixed point theorem, we deduce the existence of a unique 
solution $w(\omega, \cdot)$ of \eqref{NS2} in the unit ball  of
$H^1_{uloc}(\Omega^\eps(\omega))$, for all $\omega \in P$. Back to the
original variables, there exists a solution $u^\eps(\omega, \cdot)$
 of \eqref{NS}, which is unique in the ball of center $u^0$ and radius one in
 $H^1_{uloc}(\Omega^\eps(\omega))$.

Uniqueness in the space $B_2(\Omega^\eps(\omega))$ is deduced  easily from 
 \cite[theorem 2.3, p739]{Ladyvzenskaja:1983}. The idea is still to obtain
 local estimates on the difference of two solutions. Due to the quadratic
 term in Navier-Stokes equations, the  
inequality \eqref{Gronwallunique}
 is modified by a nonlinear term, but still leads
to uniqueness.  We do not give further details, and refer to 
 \cite{Ladyvzenskaja:1983}.

%%%%%%%%%%%%%%%%%%%%%%%%%%%%%%%%%%%%%%%%%%%%%%%%%%%%%
\subsubsection{Measurability} \label{secmeasure}

To conclude the proof of theorem \ref{existenceNS}, it remains to  check
the measurability properties of $\omega \mapsto \tilde{u}^\eps(\omega,
\cdot)$,  $P \mapsto H^1_{loc}(\R^2)$,  where the $\quad \tilde{} \quad$ stands
for the extension by zero outside $\Omega^\eps$. We know from the previous
section  that $w = u^\eps - u^0$  is the fixed point of a contraction.
 Precisely, for all $\omega \in P$, $w(\omega, \cdot)$  
 is the strong limit in
$H^1_{uloc}(\Omega^\eps(\omega))$   of the sequence $(w^k(\omega,
\cdot))_{k \in \N}$, satisfying   
\begin{equation} \label{NS3}
\left\{
\begin{aligned} 
& w^k \cdot \na w^{k+1} + u^0 \cdot \na w^{k+1} + w^{k+1} \cdot \na
 u^0 + \na q^{k+1} - \nu  \Delta w^{k+1}
  = f^\eps,  \\
&  x \in \Omega^\eps\setminus\left(\Sigma_0 \cup \Sigma_1\right), \quad 
 \div w^{k+1} = 0, \quad  x \in \Omega^\eps, \\
&  [w^{k+1}]\vert_{\Sigma_0 \cup \Sigma_1}
  =0, \quad \left[ \pa_2 w - q e_2 \right]\vert_{\Sigma_0
    \cup \Sigma_1} = -
 6\phi\\
&  \int_{\sigma(x_1)} \!\!\!\! w^{k+1}_1 \:= \: 0, \quad  
w^{k+1}\vert_{\pa \Omega^\eps} = 0,
\end{aligned}
\right.
\end{equation}
in which  we take 
$ \|w^0 \|_{H^1_{uloc}} \le 1$, with $\omega \mapsto \tilde{w}^0(\omega,
\cdot)$   measurable. 
Thus, it is enough to show 
that for all $k$, $\omega \mapsto\tilde{w}^k(\omega, \cdot)$
 is measurable from $P$ to $H^1_{loc}(\R^2)$.
Therefore, we will prove that for any function 
$$v(\omega, \cdot) \in 
H^1_{uloc}(\Omega^\eps(\omega)), \quad   \omega \mapsto 
\tilde{v}(\omega, \cdot) \mbox{ measurable, }$$
the solution $w$ of \eqref{linearNS} (with $f = f^\eps$,
 $G=0$, $\tilde{\phi} = -6\phi$) is such that $\omega \mapsto
 \tilde{w}(\omega, \cdot)$ is measurable. 

With  previous arguments and notations, for all $\omega$,  there is a
 subsequence $\displaystyle 
w_{\phi_\omega(n)}(\omega, \cdot)$ of $w_n(\omega, \cdot)$
 that converges weakly in  $\displaystyle H^1_{loc}(\omega^\eps(\omega))$ 
to $w(\omega, \cdot)$.  
By uniqueness of $w$ as a
solution of \eqref{linearNS2}, one can easily check  that the whole
sequel $w_n(\omega, \cdot)$ converges to $w(\omega, \cdot)$. Hence, by
Pettis theorem \cite{Yosida:1995}, it is enough to show the measurability of $\omega \mapsto
\tilde{w}_n(\omega, \cdot)$.
   
 To do so, we follow ideas of \cite{Abddaimi:1996} on  an
 elliptic problem from homogeneization. Let  us introduce 
$$C_n  = ]-n,n[ \times ]-2,2[,$$ 
and the spaces 
\begin{align*}
V_n & = \left\{ w \in H^1_0\left(C_n\right), \: \div w = 0 
    \right\}, \\
V_n(\omega)  & =  \left\{ w \in V_n, \: w = 0 \mbox{ in }
    C_n\setminus\Omega^\eps(\omega,n) \right\}. 
\end{align*}
Let $\pi_n(\omega) : V_n \mapsto V_n$  the orthogonal projection on 
 $V_n(\omega)$. By definition of $(P, {\cal P})$, it is easy to show that 
{\em for all $n$, the set-valued map 
\begin{equation} \label{setvalued}
{\cal F}_n : \:  (P, {\cal P}) \mapsto {\cal P}_c(H^1(\R^2)), \quad
 \omega \mapsto  V_n(\omega), 
\end{equation}
is measurable.} Following \cite{Castaing:1977}, 
we remind that a set-valued map
$$ {\cal F } : (P, {\cal P}) \mapsto {\cal P}_c(X), \quad \omega \mapsto
{\cal F}(\omega), $$
from a measurable space $P$ to  the non-empty complete subsets of a separable 
metric space $X$ is measurable if: for all open subset ${\cal O}$ of $X$, 
the set 
$$ \left\{ \omega, {\cal F}(\omega) \cap {\cal O} \neq \emptyset
\right\} $$
belongs to ${\cal P}$ (in our case, this set is open, and so belongs to
the $\sigma-$algebra ${\cal P}$). 
 For all $w$ in $V_n$, there exists by property
\eqref{setvalued}  a sequence of measurable selections 
$\sigma_j(\omega)$ such that $
\sigma_j(\omega) \rightarrow \pi_n(\omega) w$ strongly in $V_n$ ({\it c.f. 
\cite[theorem III.9, p67]{Castaing:1977}}). We thus
get, for all $w$, $w'$ in $V_n$, 
$$ (( \pi_n(\omega) u, v )) \: = \: \lim_j (( \sigma_j(\omega), w')),$$ 
where $(( \cdot, \cdot ))$ denotes the usual scalar product on
$H^1_0(C_n)$. Thus, for all $w$, $w'$, 
$$\omega \mapsto (( \pi_n(\omega) w, \, w' ))$$
 is measurable, which by Pettis theorem yields the measurability
 of $\pi_n$. 

Assume now  that $\omega \mapsto \tilde{v}(\omega, \cdot)$ is measurable.
 There exists $\alpha_n > 0 $, such that for all $\omega$, 
the solution $w_n(\omega, \cdot)$ of \eqref{linearNS} is the fixed point of the
 contraction: 
$$ F(\omega): V_n \mapsto V_n, \quad  w \mapsto \pi_n(\omega) \biggl( 
\alpha_n \Bigl(  A(\omega)w - l(\omega) \Bigr) + w \biggr), $$ 
where $l(\omega) \in V_n$, and $A(\omega) : V_n \mapsto V_n$  are defined by 
\begin{align*}
 (( \, l(\omega), w' \, )) = & 6\phi \int_{\Sigma_0(n)} w' - 
6\phi \int_{\Sigma_1(n)} w' + \int_{C_n} \tilde{f}^\eps(\omega, \cdot) \, w', \\
 ((\, A(\omega)w, w' \,)) = & \nu \int_{C_n} \Bigl( \na w \cdot \na w' +
(\tilde{v}(\omega, \cdot) + \tilde{u}^0(\omega, \cdot) )
 \cdot \na w \cdot w' \\ 
& + w \cdot \na u^0(\omega, \cdot) \cdot w' \Bigr) 
\end{align*}
Thus, $w_n(\omega, \cdot)$
 is obtained as the limit of the sequence 
$$w_n^{j+1} = F(\omega)(w_n^j).$$ 
It is clear that $\omega \mapsto l(\omega)$ and $\omega \mapsto A(\omega)$
are measurable. Together, with the measurability of $\pi_n(\omega)$, it
follows easily  that $\omega \mapsto F(\omega) w$ is measurable
 for all $w$ in $H^1(\R^2)$, hence the application $\omega \mapsto
 \tilde{w}^j_n(\omega, \cdot)$, from $P$ to $H^1(\R^2)$, for all $j$. As
 $\tilde{w}_n$ is the limit of $\tilde{w}^j_n$, theorem \ref{existenceNS}
 is proved.

%%%%%%%%%%%%%%%%%%%%%%%%%%%%%%%%%%%%%%%%%%%%%%%%%%%%%%%%%%%%%%%%%%%%%%%%%%%%
\subsection{Estimates for Dirichlet wall law}

We now turn to the proof of theorem \ref{ThDirichlet}. Instead of
estimates in $B_2$ norm, we will show the following more precise
$L^2_{uloc}$ estimates, valid for $\phi$ small enough:
\begin{align} 
\label{Dirichlet_uloc1}
& \|\nabla w\|_{L^2_{uloc}(\Omega^\eps)}  \leq C\sqrt{\eps} \\
\label{Dirichlet_uloc2}
& \|w\|_{L^2_{uloc}(\Sigma_{0,1})}  \leq C \eps \\
\label{Dirichlet_uloc3}
& \|w\|_{L^2_{uloc}(\Omega)}  \leq C \eps.
\end{align}
In fact,  the first of these inequalities  is just
\eqref{h1uloc} with $G=0$, and the second one
 is an easy consequence of \eqref{Poincare2}. 
So we just need to focus
on the third inequality, which will be proved through  a duality argument.
The main task is to estimate $\|w\|_{L^2(\Omega_R)}$ (where
$\Omega_R= \Omega \cap \{-R<x_1<R\}$). In $\Omega$, we can see $w$ as the
solution of a modified Stokes system, with inhomogeneous Dirichlet
conditions on $\Sigma_0$ and $\Sigma_1$ :
\begin{equation}
\left\{\begin{aligned}
&\nabla q -\nu \Delta w =  -(w+u^0)\cdot\nabla w - w\cdot\nabla u^0, 
    \quad x\in \Omega \\
& \div w= 0, \quad x\in \Omega \\
& w \text{ given on $\Sigma_{0,1}$,}
\end{aligned} \right.
\end{equation}
the boundary values of $w$ being controlled by \eqref{Dirichlet_uloc2}.
We introduce the adjoint equations in $\Omega_R$, that is the
Stokes system with vanishing boundary conditions and  non-zero source term:
\begin{equation} \label{Stokes} 
\left\{ \begin{aligned}
& \nabla \pi -\nu \Delta v = \varphi, \quad x\in\Omega_R \\
& \div v=0, \quad x\in\Omega_R, \\
& v|_{\partial\Omega_R}=0,
\end{aligned} \right.
\end{equation}
where $\varphi \in L^2(\Omega_R)$. We proceed in two steps. First, we 
establish two regularity lemmas for  \eqref{Stokes}. Then, by
appropriate choice of $\varphi$,  we obtain the claimed 
 estimate \eqref{Dirichlet_uloc3}.

%%%%%%%%%%%%%%%%%%%%%%%%%%%%%%%%%%%%%%%%%%%%%%%%%%%%%%%%
\subsubsection{Two lemmas for the Stokes system}

It is well known from the elliptic regularity theory 
that the solution $(v,\pi)$ of \eqref{Stokes} lies in $H^2
\times H^1$ as soon as $\varphi\in L^2$, and that there exists
$C(\Omega_R)$ such that $\|v\|_{H^2(\Omega_R)} + \|\pi\|_{H^1(\Omega_R)} \leq
C(\Omega_R) \|\varphi\|_{L^2(\Omega_R)}$ (here $\pi$ is the zero-mean
determination of the pressure).
However, the constant $C(\Omega_R)$ depends {\it a priori} on $R$.
 We need an estimate independent of $R$, 
which may be obtained by localization.

\begin{lemma} \label{lem_Stokes1}
There exists $C>0$ such that for all $R>1$, the solution $(v,\pi)$ of
\eqref{Stokes} satisfies
$$\|v\|_{H^2(\Omega_R)} + \|\nabla \pi\|_{L^2(\Omega_R)} \leq C
\|\varphi\|_{L^2(\Omega_R)}.$$
\end{lemma}
\textit{Proof.}
 
Let $\dis \chi \in C^\infty(\R)$ such that $\chi=1$ on $[0,1]$, $\chi=0$ outside
$[-1,2]$. Let $\chi_k(x_1)=\chi(x_1-k)$. We will also make use of two other
truncation functions $\theta_L$ et $\theta_R$ such that $\dis 
\theta_L \in C^\infty(\R)$, $\theta_L=1$ on $]-\infty,2]$, $\theta_L=0$ on $[3,\infty[$,
and $\theta_R(x_1)=\theta_L(-x_1)$. Let $K=\lfloor R \rfloor -1$.
We localize the equations thanks to the functions $\theta_L$, $\chi_k$ and
$\theta_R$: let $v_k(x_1,x_2)=\chi_k(x_1) v(x_1,x_2)$ ($-K\leq k \leq
K-1$), $v_L=\theta_L(x_1+R) v$, $v_R=\theta_R(x_1-R) v$. The function $v_k$
(resp. $v_L$, $v_R$) vanishes outside the domain $\dis 
D_k=(k-1,k+2)\times(0,1)$
(resp. $\dis D_L=(-R,-R+3)\times(0,1)$, $\dis 
D_R=(R-3,R)\times(0,1)$). Let $c_k$ be 
the mean of the pressure $\pi$ on $D_k$. We set $\pi_k=\chi_k
(\pi-c_k)$. We define similarly $\pi_L$, $\pi_R$.
The functions $v_k,\pi_k$ solve the following equations in $D_k$ :
\begin{equation}
\left\{ \begin{aligned}
& \nabla \pi_k - \nu \Delta v_k = \varphi_k,   \quad x\in D_k \\
& \div v_k = \gamma_k,   \quad x\in D_k \\
& v_k|_{\pa D_k}=0,
\end{aligned} \right.
\end{equation}
where 
\begin{gather*}
\varphi_k=\chi_k \varphi - \nu(2\nabla\chi_k\cdot\nabla v+ v\cdot\Delta
\chi_k) +\nabla \chi_k (\pi-c_k), \\
\gamma_k=\nabla \chi_k \cdot v.
\end{gather*}
In $D_L$ and $D_R$, $v_L$ and $v_R$ solve  analogous systems (note
that $v_L$  equals zero on the left boundary of $D_L$ , for $v$ does).
%The energy estimate for this system reads
%$$\|\nabla v_k\|_{L^2(D_k)}^2 \leq C(\|\varphi_k\|_{L^2(D_k)}^2 +
%\|\gamma_k\|_{L^2(D_k)}^2),$$
The classical elliptic theory applied to this system yields
$$\|v_k\|_{H^2(D_k)}^2 + \|\nabla \pi_k\|_{L^2(D_k)}^2 \leq
C(\|\varphi_k\|_{L^2(D_k)}^2 + \|\gamma_k\|_{H^1(D_k)}^2).$$
Of course, the constant $C$ is independent of $k$.
We have
\begin{align*}
\sum_{k=-K}^{K-1} \|\varphi_k\|_{L^2(D_k)}^2 + \|\gamma_k\|_{H^1(D_k)}^2
& \leq C(\|\varphi\|_{L^2(\Omega_R)}^2 + \|v\|_{H^1(\Omega_R)}^2) \\
& \quad {} + C \sum_k \|\pi-c_k\|_{L^2(D_k)}^2.
\end{align*}
We can control $\pi-c_k$ in the following way. One can find $C>0$ such
that for all  $f\in L^2(D_k)$ with $\int_{D_k} f=0$, there exists a
function $\psi \in H^1_0(D_k)$ satisfying $\div \psi =f$ and
    $\|\psi\|_{H^1_0(D_k)} \leq C\|f\|_{L^2(D_k)}$. Hence we have for such
$f$ and $\psi$,
\begin{align*}
\biggl|\int_{D_k} (\pi-c_k) f \biggr| & = \biggl| \int_{D_k} \nabla\pi 
    \cdot \psi \biggr|  \\
& \leq \|\nabla\pi\|_{H^{-1}(D_k)} \|\psi\|_{H^1_0(D_k)}  \\
& \leq C\|\na \pi\|_{H^{-1}(D_k)} \|f\|_{L^2(D_k)}.
\end{align*}
Thus we obtain
$$\|\pi-c_k\|_{L^2(D_k)} \leq C\|\nabla \pi\|_{H^{-1}(D_k)} \leq
    C(\|\varphi\|_{L^2(D_k)} + \|\nabla v\|_{L^2(D_k)}).$$
We have similar results for $v_L$ and $v_R$. 
Noticing that on $(k,k+1)\times (0,1)$, we have $v_k=v$ and $\nabla
\pi_k=\nabla p$, and similarly with $v_L$, $v_R$ on
$(-R,-R+2)\times(0,1)$ and $(R-2,R)\times(0,1)$, the previous
inequalities lead  to 
\begin{align*}
\|v\|_{H^2(\Omega_R)}^2 + \|\nabla \pi\|_{L^2(\Omega_R)}^2  & \leq
     \sum_{k=-K-1}^K 
 (\|v_k\|_{H^2(D_k)}^2 + \|\nabla \pi_k\|_{L^2(D_k)}^2) \\
&  \leq  \:  C (\|\varphi\|_{L^2(\Omega_R)}^2 + \|v\|_{H^1(\Omega_R)}^2) \\
&  \leq  \:  C \|\varphi\|_{L^2(\Omega_R)}^2
\end{align*}
(where indices $k=-K-1$ and $k=K$ stand for  $k=L$ and $k=K$ to
lighten  notations). Last bound comes from Poincar\'e inequality and 
 from the standard  energy estimate of \eqref{Stokes}
in the whole domain $\Omega_R$. 
This ends the proof of the lemma.

\bigskip

Our second lemma is an $H^1_{uloc}$ estimate for the Stokes system in
$\Omega_R$, uniform in $R$. We set 
$$\|v\|_{L^2_{uloc}(\Omega_R)}^2=\sup_{-R<\tau<R-1} \int_{(\tau,\tau+1)
    \times(0,1)} |v|^2$$
and $\|v\|_{H^1_{uloc}(\Omega_R)}^2 = \|v\|_{L^2_{uloc}(\Omega_R)}^2 + 
\|\nabla v\|_{L^2_{uloc}(\Omega_R)}^2$.

\begin{lemma}
There exists a constant $C>0$ such that for all $R>1$, we have
$$\|v\|_{H^1_{uloc}(\Omega_R)} \leq C \|\varphi\|_{L^2_{uloc}(\Omega_R)}.$$
\end{lemma}
\textit{Proof.} 

The proof is very similar to the calculations which led to the estimate
\eqref{h1uloc}, so we only sketch it briefly.
We write an energy estimate on the domain $(-R,r)\times(0,1)$, and integrate
with respect to $r$, from $\eta$ to $\eta+1$, to obtain a backward Gronwall
type inequality:
$$F(\eta) \leq C\int_{(-R,\eta+1)\times(0,1)} |\varphi|^2 + C F'(\eta) \leq
     C\|\varphi\|_{L^2_{uloc}}^2 (R+\eta+1) + C F'(\eta)$$
($-R \leq \eta \leq R-1$), where 
$$F(\eta) = \int_\eta^{\eta+1} \int_{(-R,r)\times(0,1)} |\nabla v|^2 \,dx
     \,dr.$$
For $\eta=R-1$, we have 
$$F(R-1) \leq \|\nabla v\|_{L^2(\Omega_R)}^2 \leq
     \|\varphi\|_{L^2(\Omega_R)}^2 \leq 2R\|\varphi\|_{L^2_{uloc}}^2,$$ 
thus we obtain
$$F(\eta) \leq C\|\varphi\|_{L^2_{uloc}}^2(R+\eta+1).$$
Similarly, we have, for $-R \leq \eta \leq R-1$,
$$G(\eta) = \int_{\eta}^{\eta+1} \int_{(r,R)\times(0,1)} |\nabla v|^2 \,dx 
     \,dr \leq C\|\varphi\|_{L^2_{uloc}}^2(R-\eta).$$
Now let $\tau \in [-R+1,R-1]$. Again, an energy estimate on 
$(\tau-r,\tau+r) \times (0,1)$ leads to
$$H_\tau(\eta) = \int_{\eta}^{\eta+1} \int_{(\tau-r,\tau+r) \times (0,1)}
     |\na v|^2 \,dx \,dr \leq C\|\varphi\|_{L^2_{uloc}}^2(\eta+1) + C
     H_\tau'(\eta)$$ 
($0 \leq \eta \leq R-1-|\tau|$). If, for instance, $\tau>0$, we have
\begin{align*}
H_\tau(R-1-\tau) & \leq \int_{(2\tau-R,R) \times (0,1)} |\na v|^2 \\
& \leq  G(2\tau-R-1) \\
& \leq C\|\varphi\|_{L^2_{uloc}}^2(2R-2\tau+1),
\end{align*}
hence we conclude that for all $0\leq \eta \leq R-1-\tau$,
$$H_\tau(\eta) \leq C\|\varphi\|_{L^2_{uloc}}^2 (\eta+1).$$
If $\tau<0$, we obtain a similar conclusion using $F$ instead of $G$.
For $\eta=1/2$, we get 
$$\int_{(\tau-\frac{1}{2},\tau+\frac{1}{2})\times (0,1)}\!\!\!\!\!\!\!\!\!
 |\nabla v|^2 \leq
     C \|\varphi\|_{L^2_{uloc}}^2.$$ 
The result follows.

%%%%%%%%%%%%%%%%%%%%%%%%%%%%%%%%%%%%%%%%%%%%%%%%%%%%%%%
\subsubsection{The duality estimate}

Now we estimate $w$ using a duality method. 
Let $R>0$, and let $v^R$, $\pi^R$ denote the solution of the Stokes problem
\eqref{Stokes} in $\Omega_R$ with source term $\varphi=w|_{\Omega_R}$.
Using these auxiliary functions, let us show the following inequality for
$\phi$ small enough: 
\begin{equation} \label{duality}
\int_{\Omega_R} |w|^2 \leq CR\eps^2 + C\int_{|x_1|=R} |w|^2 \,dx_2.
\end{equation}
To prove this, our starting point will be:
\begin{align} \label{dual1}
\int_{\Omega_R} |w|^2 & = \int_{\Omega_R} w\cdot(-\nu \Delta v^R 
       + \nabla \pi^R) \notag \\
& = \int_{\Omega_R} (-\nu \Delta w \cdot v^R + w\cdot \na\pi^R)
       - \int_{\pa \Omega_R} \nu w\cdot \frac{\pa v^R}{\pa n}.
\end{align}
We have to control the various terms appearing in \eqref{dual1}.
We recall the equation satisfied by $w$ in $\Omega_R$:
$$-\nu\Delta w + (w\cdot\na)(w+u^0) + u^0_1\pa_1 w + \na q = 0.$$
Thus integrations by parts lead to
\begin{align} \label{dual2}
\left| \int_{\Omega_R} \nu \Delta w \cdot v^R \right| & = \left| 
     \int_{\Omega_R} -w\cdot (w\cdot\na)v^R +v^R\cdot(w\cdot\na)u^0 -u^0_1
     w\cdot\pa_1 v^R \right| \notag \\
& \leq \int_{\Omega_R} |w|^2 |\nabla v^R| + C\phi \|w\|_{L^2(\Omega_R)} 
     \|v^R\|_{H^1(\Omega_R)}
\end{align}
Let $(a_0,...,a_N)$ be a regular subdivision of the interval $[-R,R]$, such
that $a_0=-R$, $a_N=R$, and $1\leq a_{k+1}-a_k \leq 2$ (hence $N\leq 2R$). 
We divide the
domain $\Omega_R$ into the corresponding cells $C_k=(a_{k-1},a_k) \times 
(0,1)$. Then we can write, using the Gagliardo-Nirenberg inequality in 
each cell $C_k$,
\begin{align} \label{dual3}
\int_{\Omega_R} |w|^2 |\nabla v^R| & \leq \sum_{k=1}^N \|w\|_{L^4(C_k)}^2
     \|\nabla v^R\|_{L^2(C_k)} \notag\\
& \leq \sum_k C\|w\|_{L^2(C_k)} \|w\|_{H^1(C_k)} \|\nabla
     v^R\|_{L^2_{uloc}} \notag\\ 
& \leq \sum_k (\eta\|w\|_{L^2(C_k)}^2 + C(4\eta)^{-1} \|w\|_{H^1_{uloc}}^2
     \|v^R\|_{H^1_{uloc}}^2) \notag\\
& \leq \eta \|w\|_{L^2(\Omega_R)}^2 + CR\eta^{-1} \|w\|_{H^1_{uloc}}^4,
\end{align}
thanks to the second lemma, with $\eta$ a small parameter to be chosen
later. We obtain, using the energy inequality for 
\eqref{Stokes}, estimate \eqref{Dirichlet_uloc1} and
the previous inequalities \eqref{dual2}, \eqref{dual3}:
\begin{equation} \label{dual4}
\left| \int_{\Omega_R} \nu \Delta w \cdot v^R \right| \leq (\eta+C\phi)
\|w\|_{L^2(\Omega_R)}^2 + CR\eta^{-1} \eps^2.
\end{equation}
Next, we treat the boundary integral in \eqref{dual1}.
We have, for $i=0$, 1, thanks to \eqref{Dirichlet_uloc2} and lemma
\ref{lem_Stokes1},
\begin{align} \label{dual5}
\left| \int_{\Sigma_i(R)} \nu w\cdot \pa_2 v^R \,dx_1 \right| & \leq \nu 
       \|w\|_{L^2(\Sigma_i(R))} \|\nabla v^R\|_{L^2(\Sigma_i(R))} \notag\\
& \leq C \eps \sqrt{R} \|v^R\|_{H^2(\Omega_R)} \notag\\
& \leq C\eta^{-1} \eps^2 R + \eta \|w\|_{L^2(\Omega_R)}^2.
\end{align}
We also write, using again a trace theorem,
\begin{align} \label{dual6}
\left|\int_{|x_1|=R} \nu w\cdot \pa_1 v^R \,dx_2 \right| & \leq C \left( 
    \int_{|x_1|=R} |w|^2 \,dx_2\right)^{\frac{1}{2}}
    \|v^R\|_{H^1(\Omega_{R-1,R})} \notag \notag\\
& \leq C\eta^{-1} \int_{|x_1|=R} |w|^2 \,dx_2 + \eta \|w\|_{L^2(\Omega_R)}^2
\end{align}
(here the integrals on the sections $x_1=\pm R$ are restricted to the
portion of each section contained in $\Omega$).
It remains to deal with the integral involving the pressure $\pi^R$ in
\eqref{dual1}. 
We use again the subdivision of $\Omega_R$ into the cells $C_k$. Let $c_k$
denote the mean of $\pi^R$ on $C_k$. We integrate by parts $w\cdot\na\pi^R
=w\cdot\na(\pi^R-c_k)$ on each cell $C_k$:
\begin{align} \label{dual7}
\int_{\Omega_R} w\cdot\na\pi^R  & = \sum_{k} \int_{\pa C_k} w\cdot n
      (\pi^R-c_k) \notag\\
& = \sum_{k} \biggl( \int_{\Sigma_1 \cap \pa C_k}  w_2(\pi^R-c_k) \,dx_1
      - \int_{\Sigma_0 \cap \pa C_k} w_2 (\pi^R -c_k) \,dx_1 \biggr)
      \notag\\ 
& \quad\ {} - \int_{x_1=-R} w_1 (\pi^R-c_1) \,dx_2 + \int_{x_1=R}
      w_1(\pi^R-c_N) \,dx_2 \notag\\ 
& \quad\ {} + \sum_{k=1}^{n-1} \int_{x_1=a_k} w_1(c_{k+1} - c_k)
      \,dx_2. 
\end{align}
We bound the quantities appearing in the r.h.s. of \eqref{dual7}.
We have, for $i=0$, 1,
\begin{align*}
\sum_{k} \left|\int_{\Sigma_i \cap \pa C_k} w_2(\pi^R-c_k) \,dx_1 \right|
    & \leq \sum_k C \|w\|_{L^2_{uloc}(\Sigma_i)} \|\pi^R-c_k\|_{H^1(C_k)} \\
& \leq \sum_k (C\eta^{-1} \eps^2 + \eta \|\nabla\pi^R\|_{L^2(C_k)}) \\
& \leq C\eta^{-1} \eps^2 R + \eta \|w\|_{L^2(\Omega_R)}^2.
\end{align*}
using \eqref{Dirichlet_uloc2}, the Poincar\'e-Wirtinger inequality applied
to the zero-mean function $\pi^R-c_k$ in $C_k$, and lemma \ref{lem_Stokes1}.
Similarly, we have for the two terms corresponding to $x_1=\pm R$:
$$\left| \int_{x_1=\pm R} w^1(\pi^R - c_{1,N}) \,dx_2 \right| \leq
      C\eta^{-1} \int_{|x_1|=R} |w|^2 \,dx_2 +
      \eta\|w\|_{L^2(\Omega_R)}^2.$$  
To deal with the last term of \eqref{dual7}, we use the fact that
$\int_{\sigma(a_k)} w_1 \,dx_2=0$, where $\sigma(a_k)=\{(a_k,x_2) \in
\Omega^\eps \}$ (here we take into account 
the exterior $\Omega^\eps\setminus\Omega$ of $\Omega$). Hence
\begin{align*}
\left|\int_{x_1=a_k} w_1 \,dx_2 \right| & = \left| \int_{\sigma(a_k)
      \setminus \Omega} w_1 \,dx_2 \right| \\
& \leq \left(\eps \int_{\sigma(a_k) \setminus \Omega} |w_1|^2 \,dx_2
      \right)^{1/2} \\
& \leq C \sqrt{\eps} \|w\|_{H^1_{uloc}(\Omega^\eps)}
\end{align*}
(we can apply the trace theorem in $\Omega^\eps \cap \{a_k<x_1<a_{k+1}\}$
for instance, with a constant independent of $k$ since the boundaries of
$\Omega^\eps$ are $K$-Lipschitz continuous). Moreover, we obviously have
$|c_{k+1}-c_k|\leq C\|\nabla\pi^R \|_{L^2(C_k \cup C_{k+1})}$, thus
estimate \eqref{Dirichlet_uloc1} and lemma \ref{lem_Stokes1} again imply
\begin{align*}
\left|\sum \int_{x_1=a_k} w_1(c_{k+1}-c_k) \,dx_2 \right| & \leq  \sum
     C\eps \|\nabla\pi^R \|_{L^2(C_k \cup C_{k+1})} \\
& \leq C\eta^{-1} \eps^2 R + \eta \|w\|_{L^2(\Omega_R)}^2.
\end{align*}
Up to this point, the r.h.s. of \eqref{dual7} can be controlled as follows:
\begin{equation} \label{dual8}
\left| \int_{\Omega_R} w\cdot\nabla \pi^R \right| \leq C\eta^{-1} \eps^2 R
     + 4\eta \|w\|_{L^2(\Omega_R)}^2 + C \int_{|x_1|=R} |w|^2 \,dx_2.
\end{equation}
Then inequalities \eqref{dual4}, \eqref{dual5}, \eqref{dual6} and
\eqref{dual8} enable us to write, from \eqref{dual1}:
$$\int_{\Omega_R} |w|^2 \leq (C\phi+8\eta)\|w\|_{L^2(\Omega_R)}^2 + 
         CR\eta^{-1} \eps^2 + C\eta^{-1} \int_{|x_1|=R} |w|^2 \,dx_2.$$
We take $\eta=1/10$ and $\phi<\phi_0$ such that $C\phi_0<1/10$. This leads
to \eqref{duality}.

Let $z(R)=\|w\|_{L^2(\Omega_R)}^2$. Inequality \eqref{duality} can be
rewritten $z(R) \leq {\cal C} R \eps^2 + {\cal C} z'(R)$. Following
\cite{Ladyvzenskaja:1983}, we deduce the  estimate
\eqref{estimatesDirichlet} from this differential inequality and from the
fact that $z(t)$ has subexponential growth at infinity (Poincar\'e
inequality and \eqref{Dirichlet_uloc1} imply that it has linear growth).
Indeed, let $\dis \tilde z(R)={\cal C} R\eps^2 + {\cal C}^2 \eps^2$;
we have  $
z'(R)-\tilde z'(R) \geq \frac{z(R)-\tilde z(R)}{\cal C}$. Hence if 
there existed $R_0$ such that $\dis z(R_0) -\tilde z(R_0)=A>0$, we would have
$z(R)-\tilde z(R) \geq Ae^{(R-R_0)/{\cal C}}$ for all $R \geq R_0$, and $z$
would 
grow exponentially fast at infinity, which is not the case. Thus we obtain
for all $R>1$, $\dis z(R) \leq {\cal C} R\eps^2 + {\cal C}^2 \eps^2$, that is
\eqref{estimatesDirichlet}. 

Finally, we show \eqref{Dirichlet_uloc3}. It is a direct consequence of
what we have already proved. For $R=1$, we have $\|w\|_{L^2(\Omega_1)}^2 =
z(1) \leq C\eps^2$. Now we can repeat the whole proof but with another
origin: take $x_1=\tau$ instead of $x_1=0$ as origin, replace $\Omega_R$ with
$\Omega(\tau,R) =\Omega_R + (\tau,0)$. This yields $z_\tau(R)=
\|w\|_{L^2(\Omega(\tau,R))}^2 \leq {\cal C} R\eps^2 + {\cal C}^2 \eps^2$,
\emph{with the same constant $\cal C$ as before}. Hence we have again 
$\|w\|_{L^2(\Omega(\tau,1))}^2 \leq C\eps^2$, and we conclude that
$\|w\|_{L^2_{uloc}(\Omega)}^2 \leq C \eps^2$.

%%%%%%%%%%%%%%%%%%%%%%%%%%%%%%%%%%%%%%%%%%%%%%%%%%%%%%%%%%%%%%%%%%%%%%%%%%%%
%%%%%%%%%%%%%%%%%%%%%%%%%%%%%%%%%%%%%%%%%%%%%%%%%%%%%%%%%%%%%%%%%%%%%%%%%%%% 
\section{The boundary layer analysis} \label{sectionBL}
\subsection{Formal expansion}
To obtain a refined wall law at $\pa \Omega$, one must clarify  
the structure of the flow near the rough boundary. Namely,  
we will show that $u^\eps$ behaves as 
  \begin{equation} \label{Ansatz2}  
\begin{aligned}
u^\eps(\omega, x) \approx  u^\eps_{app}(\omega,x) & = u^0(\omega, x)  +
\eps  U_l\left(\omega,  
\frac{x_1}{\eps}, \frac{x_2}{\eps} \right) + \eps 
 U_u\left(\omega, \frac{x_1}{\eps},
\frac{x_2 - 1}{\eps} \right)  \\
&\quad + \eps  u^1(\omega,x).  
\end{aligned}  
\end{equation}
In this formal expansion,  $U_{l,u}(\omega,y)$ are {\em boundary layer
profiles}, expressing the strong gradients of $u^\eps$ near the
boundary. They are
defined on the rescaled domains ${\cal R}_{l,u}(\omega)$. They should be
localized near the boundary, and cancel  the jump of the normal derivative
 of $u^0$. As usual, to derive  the system they satisfy, we plug expansions
 \eqref{Ansatz2}  in 
 equations \eqref{NS}, and collect terms of order $\eps^{-1}$. Formally,
 this  leads to  the Stokes system  \eqref{BL}a,b. Jump conditions 
\eqref{BL}c,d ensure the regularity of the approximation through $\pa
\Omega$.  Finally, \eqref{BL}d is in agreement with the 
 Dirichlet condition \eqref{BC}. 

 As mentioned in theorem \ref{existenceBL},  we will show that 
$U_{l,u}(\omega, \cdot)$ converge to some constant fields 
$$U_{l,u}^\infty(\omega) = (U^\infty_{l,u,1}(\omega), 0)$$
as $y_2 \rightarrow +\infty$. Because of these fields, which do not match
the required flux and boundary conditions, there is an  additional 
corrector $u^1(\omega, x)$. It is given by 
\begin{equation} \label{u1} 
\left\{
\begin{aligned}
u^1(\omega, x) & = \Bigl( 3 (U^\infty_{l,1} + U^\infty_{u,1}) x_2^2 -
(4 U^\infty_{l,1} +  2 U^\infty_{u,1}) x_2 - U^\infty_{u,1}, \: 0\Bigr), 
\quad x \in \Omega,\\
u^1(\omega, x) & = - U^\infty_{u}, \quad  x_2 < 0, \\
u^1(\omega, x) & = - U^\infty_{l}, \quad x_2 > 1. 
\end{aligned}
\right.
\end{equation}
Note that in the interior domain $\Omega$, $u^1$ satisfies 
\begin{equation*}
\left\{
\begin{aligned}
& u^0 \cdot \na u^1 + u^1 \cdot \na u^0 - \nu \Delta u^1 + \na p^1 = 0, \quad
x \in \Omega, \\
& \div u^1 = 0, \quad x \in \Omega, \\
& u^1\vert_{\Sigma_0} = -U^\infty_u \quad u^1\vert_{\Sigma_1} =
-U^\infty_l, \quad \int_{\sigma(x_1)}\!\!\!\!\!\!\!\!
 u^1_1 = - (U^\infty_{l,1} + U^\infty_{u,1}).
\end{aligned}
\right.
\end{equation*}

\subsection{Well-posedness of \eqref{BL}}
We now justify the previous formal computations, and prove theorem
\ref{existenceBL}. We first consider the well-posedness of system
\eqref{BL}. 
 In the study of Dirichlet wall law, the Poincar\'e
inequality \eqref{Poincare}, applied to the channel $\Omega^\eps$, 
 allowed for deterministic reasoning. In the boundary layer domains
 ${\cal R}_{l,u}(\omega)$ that are unbounded in every direction, this
 inequality does not hold, and we shall use the probabilistic modeling to
 solve \eqref{BL}. This approach borrows to homogeneization problems (see
 \cite{Jikov:1994,Beliaev:1996, Bourgeat:1994}), 
but strong changes are needed to account for the anisotropy of
 our boundary layer domains. The general idea is to construct a solution of
 the type 
\begin{equation} \label{formesolution}
U_{l,u}(\omega,y) = 
V_{l,u}\left(\tau_{y_1}(\omega), y_2 \pm h_{l,u}\circ \tau_{y_1}(\omega)\right),
\quad V_{l,u}  = V_{l,u}(\omega, \lambda),
\end{equation}
where $\tau_{y_1}$ and $h_{l,u}$ were defined in subsection
\ref{modeling}. This  requires to introduce adapted functional spaces and
 variational formulations. As the lower and upper systems are treated in
 the same way, {\em we will focus on the lower one,  
and drop  the subscript  ``l'' for brevity}. 
\subsubsection{Stochastic derivative, Convolution} 
 Let $V : P \times \R \mapsto \R^n$, $n \ge 1$, measurable. We call a
 {\em realization of V}  an  application 
$$ R_\omega[V] : \R^2 \mapsto \R^n, \quad 
y \mapsto V\left(\tau_{y_1}(\omega), y_2 + h \circ \tau_{y_1}(\omega)\right),
 \quad \omega \in P. $$ 
Notice that 
$$R_{\tau_h}(\omega)[V](y_1,y_2) \: = \: 
 R_\omega[V](y_1+h,y_2), \quad \forall h,y.$$ 
We say that $V$ is {\em  smooth} if, almost surely, $\dis R_\omega[V]$ is
smooth as a function on $\R^2$. For all smooth functions $V$, we define the
{\em stochastic derivative of $V$} as 
$$ \pa_\omega V(\omega, \lambda) =  \pa_1R_\omega[V]\left(y_1 = 0, y_2 = 
\lambda - h(\omega)\right) $$ 
We also introduce 
$$ \pa_\lambda V(\omega, \lambda)  = \pa_2R_\omega[V]\left(y_1 = 0, y_2 = \lambda
- h(\omega)\right) $$     
and the gradient $\na_\omega = (\pa_\omega, \pa_\lambda)$. Notice that,
almost surely,  
\begin{equation} \label{stochder1}
R_{\omega}\left[\na_\omega V\right] \: =  \: \na R_{\omega}[V] 
\end{equation}
which allows to define for smooth $V$ partial derivatives at any order
$\pa_\omega^\alpha \, \pa_\lambda^\beta V$. 
%Then we  introduce the space $C^\infty_c\left( P \times \R \right)$ of
%smooth functions $V$ s.t. 
%\begin{description}
%\item[i)]
%the support of $V$ is compact in $\lambda$, that is $V(\cdot,
%  \lambda) = 0$ for $|\lambda|$ large enough. 
%\item[ii)] for all  $\alpha, \beta$, 
%$  \pa_\omega^\alpha \, \pa_\lambda^\beta V \in L^\infty\left( P \times \R
%  \right).$
%\end{description}

For interesting examples of smooth functions, one can use a convolution
process. Precisely, let $\rho =
\rho(y) \in \R$ a smooth function with compact support. For all $V \in
L^1_{loc}\left(P \times \R\right) = \cap_K L^1(P \times
K)$, $K$ compact,  we define
\begin{equation} \label{convolution}
\rho * V(\omega,\lambda)  = \int_{\R^2}  \rho(y) \, 
V\left(\tau_{y_1}(\omega), \lambda
+ y_2 +   h \circ \tau_{y_1}(\omega) - h(\omega) \right) \, dy.
\end{equation}
The following lemma is easy and left to the reader. 
\begin{lemma} \label{lemmeconvolution}
 Assume  $ 1 \le p \le \infty$ and  $V \in L^p\left(P \times
  \R\right)$, resp. $L^p_{loc}$. Then $\rho * V \in L^p\left(P \times
  \R\right)$, resp. $L^p_{loc}$,  and is a smooth function. Moreover, 
the following identities hold: 
\begin{align*}
& R_\omega[\rho * V] \: = \:
\hat{\rho} * R_\omega[V], \: \mbox{ a.s. in } \omega,  \\
& \na_\omega (\rho * V) \: = \: \na \hat{\rho} * V, 
\end{align*}
where $\hat{\rho}(y) =\rho(-y)$.  Besides, if $V$ is smooth and
$\na_\omega V \in L^1_{loc}(P \times \R)$, then
$$ \na_\omega (\rho * V) = \rho * \na_\omega V. $$
\end{lemma}
%As a consequence of this lemma, for any $V \in L^\infty(P \times \R)$ with
%compact support in $\lambda$, and any $\rho$, $\rho \times V$ belongs to 
% $C^\infty_c\left( P \times \R \right)$. 
\subsubsection{Functional spaces}
We can now introduce functional spaces adapted to \eqref{BL}. Let ${\cal
  D}_0(P \times \R_+)$ the set of smooth functions $\varphi =
\varphi(\omega, \lambda) \in \R^2$,  satisfying 
\begin{description}
\item[i)] For all  $\alpha, \beta$, 
$\quad   \pa_\omega^\alpha \, \pa_\lambda^\beta \varphi \in
 L^\infty\left( P \times \R \right).$
\item[ii)] Almost surely,  the 
support of $R_\omega[\varphi]$ is included in $\{y, \: \frac{1}{\delta} > y_2 + h \circ
\tau_{y_1}(\omega) > \delta \}$ for some positive constant $\delta = \delta(\omega,
\varphi)$.  
\end{description} 
Note that realizations of functions of  ${\cal  D}_0(P \times \R_+)$
 have a support that is strictly included in ${\cal R}(\omega)$, almost
 surely in  $\omega$. 
We clearly define a scalar product on  ${\cal  D}_0(P \times \R_+)$ through 
$$(( \varphi, \tilde{\varphi} )) \: = \:  \int_{P \times \R} \na_\omega 
\varphi \cdot \na_\omega \tilde{\varphi}. $$
Let $D_0(P \times \R_+)$ the completion of ${\cal  D}_0(P \times \R_+)$ for
the norm $\| \varphi \|_{D_0} = (( \varphi,  \varphi))^{\frac{1}{2}}$. 
It is of course a Hilbert space. We state the following properties:
\begin{lemma} \label{propD0}
We have 

\smallskip
\begin{description}
\item[i)] $D_0(P \times \R_+) \: \hookrightarrow \: 
C^0(\R_+; \,  L^2(P)) \: \hookrightarrow \:  L^2_{loc}(P \times \R_+)$.
\item[ii)] For any $V$ in  $D_0(P \times \R_+)$, there exists a
  unique $D_\omega  V  \in L^2\left(P \times \R_+\right)$ such that 
\begin{equation} \label{weakder}
  \int_{P \times \R_+} \left( D_\omega V \right) \varphi
 = - \int_{P \times \R_+}  \left( \na_\omega \varphi \right) V, \quad
 \forall \varphi \in {\cal  D}_0(P \times \R_+). 
\end{equation}
Moreover,  $\: \| D_\omega V \| = \| V \|_{D_0}$. With a slight abuse of
notation, we shall write $\na_w  = (\pa_\lambda, \pa_\omega) $ instead of
$D_\omega $ (``weak stochastic derivative "). 
\item[iii)] For any $V$ in  $\dis D_0(P \times \R_+)$, almost surely,
  $R_\omega[V] \in H^1_{loc}\left({\cal R}(\omega)\right) \dis $, $
 \dis  R_\omega[V]\vert_{\pa \left({\cal R}(\omega)\right)} = 0$, $\dis \na
  R_\omega[V] = R_\omega[\na_\omega V]$. 
\end{description}
\end{lemma}
\textit{Proof.} 

{\bf i)} We can argue for $V \in {\cal D}_0$, as the
result for $D_0$ follows by density. Almost surely, $R_\omega[V]$ belongs to
$H^1_{loc}({\cal R}(\omega)$, $\quad R_\omega[V]\vert_{\pa {\cal R}(\omega)} =
0$. Thus,
\begin{equation*}
\begin{aligned}
\int_{P} |V|^2(\omega, \lambda) dP  & \: = \: \frac{1}{2} \int_{(-1,1)} dy_1 \, 
 \int_{P} |V|^2(\tau_{y_1}(\omega), \lambda) dP \\
& \: = \: \frac{1}{2} \int_{P} \, \left\| R_\omega[V]\left(y_1,\lambda -
h \circ \tau_{y_1}(\omega)\right) \right\|_{L^2(-1,1)}^2 \, dP \\
& \: \le \:
C(\lambda) \, \int_{P} \, \left\| \na R_\omega[V](x, z) \right\|_{L^2({\cal
    R}(\omega, 1))}^2 \, dP \\ 
& \: = \:  C(\lambda) \, \int_{P} \left\|  R_\omega[\na_\omega V](x, z) 
\right\|_{L^2({\cal R}(\omega, 1))}^2 \, dP \\
& \: = \: 2  C(\lambda)
 \int_{P\times \R_+} \,  |\na_\omega V|^2 dP d\lambda
\end{aligned}
\end{equation*}
where the inequality comes from classical trace theorem and Poincar\'e
inequality. 

{\bf ii)} Let $V \in D_0$, and  $V_n \in {\cal D}_0$, converging to $V$. As
$\dis \| V_n - V_m \|_{D_0} = \| \na_\omega V_n - \na_\omega V_m\|$,
 $\na_\omega V_n$ is a Cauchy sequence in $L^2(P \times \R_+)$, thus
converging to some $D_\omega V$ with $\|D_\omega V\|_{L^2} =  \| V
\|_{D_0}$. As
$$  \int_{P \times \R_+} \left( \na_\omega V_n \right) \varphi
 = - \int_{P \times \R_+}  \left( \na_\omega \varphi \right) V_n, \quad
 \forall \varphi \in {\cal  D}_0(P \times \R_+), $$
relation \eqref{weakder} follows (notice that
the right hand side converges by {\bf i)}). 

It remains to show uniqueness of $D_\omega V$. This follows easily from
the following general result:  for any $W$
in $L^p\left( P \times \R_+\right)$, resp. $L^p_{loc}$,  
 $1 \le p < \infty$, there exists
$W_n \in {\cal D}_0(P \times \R_+)$  such that 
$$ W_n \xrightarrow[n 
\rightarrow +\infty]{} W \quad \mbox{ in } L^p\left( P \times \R_+\right),
\mbox{ resp. } L^p_{loc}. $$ 
To see this,  first notice  that functions of $L^\infty(P \times \R_+)$
with compact support in $\lambda$ are dense in $L^p(P \times \R_+)$, resp.      
$L^p_{loc}$. We therefore assume that $W$ is such a function. We then
proceed by convolution: let $\rho_n$ an approximation of unity, and set 
$W_n = \rho_n * W$. From lemma \ref{lemmeconvolution}, we deduce easily
that $W_n \in {\cal D}_0(P \times \R_+)$. To obtain the convergence,
 we involve as above the realizations of $W_n$ and $W$, by extra
 integration with respect to $y_1$. Namely, for $W \in L^p$, 
\begin{align*}
 \int_{P \times \R} dP \, & \left| W_n - W \right|^p \: = \: \frac{1}{2} 
 \int_{(-1,1)} dy_1 \, \int_{P \times \R}  dP \, \left| W_n - W \right|^p \\
& \: \le  \:  \frac{1}{2}  \int_{P} \int_{(-1,1) \times \R} dy \, \left|
R_\omega\left[ \rho_n * W \right](y) -  R_\omega[W](y)  \right|^p  \\
& \: \le \: \int_P dP \, \left( \int_{(-1,1) \times \R} dy \, \left| 
\rho_n * R_\omega\left[  W \right](y) -  R_\omega[W](y)  \right|^p \right).
\end{align*}
By standard results on the (standard)
 convolution,  the integral inside the parenthesis converges to $0$ as $n
 \rightarrow +\infty$, almost surely. The dominated convergence theorem
 allows to conclude. Similar reasoning holds for $W \in L^p_{loc}$, with $P
 \times K$, $K compact$, instead of $P \times \R$. 

{\bf iii)} follows again from a density argument, and is left to the reader

\bigskip
We finally introduce the (closed) 
space of divergence-free vector fields:
$$ V_0(P \times \R_+) \: = \: \left\{ V \in D_0(P\times \R_+), \:
\na_\omega \cdot V = \pa_\omega V_1 + \pa_\lambda V_2 = 0 \right\}. $$

\subsubsection{Variational formulation}
The keypoint  to solve \eqref{BL} is to search the solution
$U(\omega,\cdot)$  as a realization: almost surely, 
\begin{equation} \label{Urealization}
 U(\omega, \cdot) = R_\omega[V], \quad V \in V_0(P \times \R_+).
\end{equation} 
At a formal level, if we substitute \eqref{Urealization} in \eqref{BL}, and
test again $R_\omega[\varphi]$,$\varphi \in V_0(P \times \R_+)$, we end
up with 
\begin{equation} \label{variational} 
\nu \int_{P \times \R_+} \na_\omega V : \na_\omega \varphi  \: = \:
 {\cal L }(\varphi),
\quad \forall \varphi \in V_0(P \times \R_+), 
\end{equation}    
where the linear form 
$$  {\cal L }(\varphi)  \: = \: 6 \phi \cdot 
\int_{P} \varphi(\omega, +h(\omega)) dP $$
comes from the inhomogeneous jump condition on the normal derivative. Note
that ${\cal L}$ is well-defined and continuous on $D_0(P \times \R_+)$ by
lemma \ref{propD0}. This suggests the following definition: 

\begin{center} \label{defvariational}
{\em $U$ is a variational solution of \eqref{BL} if it satisfies
 \eqref{Urealization}, \eqref{variational}.}
\end{center}

Notice that by Riesz theorem, there exists a unique solution $V$ to
\eqref{variational}, and so a unique variational solution $U$ to
\eqref{BL}. From lemma \ref{propD0},  one has almost surely, 
$$U(\omega, \cdot) \in  H^1\left({\cal R}(\omega)\right), \quad U(\omega,
\cdot)\vert_{\pa \left({\cal R}(\omega)\right)} = 0, \quad 
\na U(\omega,\cdot) \in L^2\left({\cal R}(\omega, R)\right), \: \forall
R, $$
with  for all $R \ge 1$, 
$$ \frac{1}{2R} \E\left(\int_{{\cal R}(\cdot, R)} |\na U(\cdot,y)|^2
dy \right) = \| \na_\omega V \|^2_{L^2} < \infty. $$
Moreover, by the ergodic theorem, we have
\begin{align*}
& \lim_{R \rightarrow +\infty} \frac{1}{R} \int_{{\cal R}(\omega,R)}  |\na
U(\cdot,y)|^2 dy  \\
 = &  \lim_{R \rightarrow +\infty}\frac{1}{R}
 \int_{(-R, R)}  dy \int_{\R_+} 
|\na_\omega V(\tau_y(\omega), \lambda)|^2 d\lambda
\end{align*}
exists almost surely, so that 
$$ \sup_{R \ge 1} \frac{1}{R} \int_{{\cal R}(\omega,R)}  |\na
U(\cdot,y)|^2 dy  < \infty \quad \mbox{almost surely }.$$

It remains to show that $U$ is almost surely a classical solution. This is
a consequence of 
\begin{lemma} 
Almost surely, $\div U = 0$ in the weak sense, and 
for every $\psi \in {\cal C}^\infty_c\left({\cal
  R}(\omega)\right)$ with  $\div \psi = 0$, we have  
$$ \nu \int_{{\cal R}(\omega)} \na U : \na \psi  = 6 \phi
 \int_{\Sigma_0} \psi(y_1,0) dy. $$
\end{lemma}
{\it Proof.} By lemma \ref{propD0}, iii),  $\div U = 0$ in the weak
sense  almost surely.  
%By definition of $V_0(P \times \R_+)$, we have 
%$$ \int_{P \times \R_+} V \cdot  \na_\omega \varphi = 0, \quad \forall 
%\varphi  \in {\cal D}_0(P \times \R_+).$$ 
%We consider again an approximation of unity $\rho_n = \rho_n(y)$, and set
%$\varphi_n = \rho_n * \varphi$. We have easily  
%\begin{align*}
%  0 & =  \int  V \cdot  \na_\omega \varphi_n =  \int V \cdot \left( \rho_n *
%\na_\omega \varphi \right) = \int \left( \hat{\rho}_n * V \right) 
% \cdot \na_\omega \varphi \\
%& = -\int  \left(\na_\omega \cdot  \left(\hat{\rho}_n * V\right) \right) \,
% \varphi. 
%\end{align*} 
%We deduce that  $\dis \na_\omega \cdot  \left(\hat{\rho}_n * V\right) =
%0$,
% and 
% applying $R_\omega$, we get 
%$$\dis \na \cdot \left( \hat{\rho}_n * R_\omega[V] \right) = 0. $$
%If we let $n$ go to infinity, we deduce that 
%$\dis  \na \cdot  R_\omega[V] = \na \cdot U = 0$ almost surely, 
%in the weak sense. 
To recover the Stokes equation in the weak sense, we reexpress the
variational equation \eqref{variational} as 
$$ 
\nu \int_{P \times \R_+} \na_\omega V : \na_\omega \varphi \: +  \:  
\int_{P \times \R_+} \na_\omega W : \na_\omega \varphi \:  \: = \:
 {\cal L }(\varphi),
\quad \forall \varphi \in D_0(P \times \R_+), 
$$ 
where $W \in V_0^\bot(P \times \R_+)$, which means that 
$$ \int_{P \times \R_+} \na_\omega W : \na_\omega \tilde{V} = 0, \quad
\forall \tilde{V} \in  V_0(P \times \R_+). $$
We consider again an approximation of unity $\rho_n = \rho_n(y)$, and set
 $\tilde{V} = \rho_n * \left(\na_\omega^\bot \varphi \right)$, 
 $\varphi \in  {\cal D}_0(P \times \R_+)$. We compute 
\begin{align*}
0 & = \int_{P \times \R_+} \na_\omega W : \na_\omega(\rho_n *
\na_\omega^\bot \varphi) 
= \int_{P \times \R_+}   \na_\omega W :  \left(\rho_n * ( \na_\omega
\na_\omega^\bot \varphi) \right) \\
& = \int_{P \times \R_+} \left( \hat{\rho}_n * (\na_\omega W) \right) : \na_\omega
\na_\omega^\bot \varphi =  - \int_{P \times \R_+} \left( \na_\omega^\bot \cdot
(\na_\omega)^2 (\hat{\rho}_n * W)\right) \, \varphi. 
\end{align*}
so that $\na_\omega^\bot \cdot
(\na_\omega)^2 (\hat{\rho}_n * W) = 0$. Applying $R_\omega$ we deduce
that 
\begin{equation} \label{curlstoch} 
  \na^\bot \cdot \Delta  \left( \rho_n * R_\omega[W] \right) = 0.
\end{equation}
Back to \eqref{variational}, we take $\varphi_n = \rho_n * \varphi$ as a test
function. One has easily the following identities:
\begin{align*}
& \int \na_\omega V : \na_\omega \varphi_n  = - \int \biggl((\na_\omega)^2
\left( \hat{\rho}_n * V \right) \biggr) \cdot \varphi \\
& \int \na_\omega W : \na_\omega \varphi_n  = - \int \biggl((\na_\omega)^2
\left( \hat{\rho}_n * W \right) \biggr) \cdot \varphi \\
& {\cal L}(\varphi_n)  = 6 \phi \cdot \int \biggl(\hat{\rho}_n *  
 \Bigl( (\omega, \lambda) \mapsto \delta(\lambda - h(\omega)) \Bigr) \biggr)
 \cdot \varphi , 
\end{align*}
where $\delta$ stands for the Dirac measure, {\it i.e.}
$$ \hat{\rho}_n * \Bigl( (\omega, \lambda) \mapsto   \delta(\lambda +
h(\omega)) \Bigr) = 
\int_\R dy_1 \, \hat{\rho}_n(y_1, -\lambda  + h(\omega)) dy_1. $$
We thus deduce from \eqref{variational} that,  
$$ - \nu (\na_\omega)^2 \left( \hat{\rho}_n * V \right) -
(\na_\omega)^2 \left( \hat{\rho}_n * W \right) = 6\phi \cdot
\hat{\rho}_n *  
 \Bigl( (\omega, \lambda) \mapsto \delta(\lambda - h(\omega)) \Bigr), $$
which implies that, almost surely: 
\begin{equation} \label{stokesregular}
-\nu \Delta \,  \left( \rho_n * U  \right)(\omega, \cdot) 
-\Delta \,  \left( \rho_n * R_\omega[W]  \right)(\omega, \cdot) = 
6 \phi \cdot \left( \rho_n *  \delta_{y_2 = 0}\right).  
\end{equation}
From \eqref{curlstoch},  we can express
$$  \Delta \,  \left( \rho_n * R_\omega[W] \right)(\omega, \cdot) = \na
p_\omega,$$
for some smooth scalar field $p_\omega$. Hence, if we multiply
\eqref{stokesregular} by a  divergence free vector field 
$\psi \in {\cal C}^\infty_c({\cal R}(\omega))$  and perform 
integration by parts,
we obtain 
 $$ \int_{{\cal R}(\omega)} \na (\rho_n *  U) : \na \psi  = \frac{\phi}{6}
 \int_{\Sigma_0} (\rho_n *\psi)(y_1,0) dy_1. $$
The limit $n \rightarrow +\infty$ gives the result. The fact that
$U(\omega. \cdot)$ is regular follows from standard ellipticity properties
of the Stokes operator. This shows the well-posedness of the boundary layer
system \eqref{BL}.

%%%%%%%%%%%%%%%%%%%%%%%%%%%%%%%%%%%%%%%%%%%%%%%%%%
\subsection{Convergence at infinity}

This section is devoted to the last part of theorem \ref{existenceBL}, that
is convergence of $U_{l,u}$ as $y_2 \rightarrow +\infty$. In the periodic
setting, the convergence is exponential, as can be seen easily from a Fourier
analysis. Namely, if the roughness has period $T$, it is straightforward
that 
$$ \left| U_{l,u}(y_1, y_2) - U_{l,u}^\infty \right| \: \le \: C
e^{-\frac{\alpha}{T}}, $$   
for some constant $\alpha$ independent of $T$. In this case, 
$U_{l,u}^\infty = \left( U_{l,u,1}^\infty, \, 0  \right)$ is just the
average of $U_{l,u}$ with respect to $y_1$. We stress that the convergence
rate $\alpha/T$ goes to zero as the period of the roughness $T$ goes to
infinity. Thus, in the random setting (in which, broadly speaking, all
periods are involved in the roughness), this analysis falls down. We will
show the convergence using ergodicity properties of $U_{l,u}$. {\em For the
  sake of brevity, we will only treat  the lower boundary layer 
  $U_l$}. Similar reasoning holds for the upper one. 

We first establish a representation formula for $U_l$, in terms of {\em 
the double layer Stokes potential.} 

\begin{lemma}
Almost surely, the solution $U_l$ of \eqref{BL} satisfies, for $y_2 > 0$, 
\begin{equation} \label{doublelayer}
U_l(\omega,y) \: = \: \int_{\R} G(y_1- y'_1,y_2) \, U_l(\omega,y_1',0) \,
dy'_1,
\end{equation}
where 
$$ G(y_1,y_2) =  \frac{2 y_2}{\pi \, \left( y_1^2 + y_2^2\right)^2} \, 
\left( \begin{matrix} y_1^2 & y_1 y_2 \\  y_1 y_2 & y_2^2 \end{matrix}
\right). $$ 
\end{lemma}
{\it Proof.}

 Let us denote by $\tilde{U}_l(\omega,y)$ the right hand side of
\eqref{doublelayer}. We remind that 
$$ U_l(\omega, y) \: = \: V(\tau_{y_1}(\omega), y_2 + h_l\circ
\tau_{y_1}(\omega)) = {\cal V}(\tau_{y_1}(\omega), y_2). $$
As a consequence,
$$\tilde{U}_l(\omega,y) = \tilde{{\cal V}}(\tau_{y_1}(\omega), y_2), $$
where 
$$ \tilde{{\cal V}}(\omega,y_2) = \int_{\R} G(y'_1,y_2)
{\cal V}(\tau_{-y'_1}(\omega),0) \, dy'_1. $$
Thus, $U_l$ and $\tilde{U}_l$ are both stationary with respect to
$y_1$. Moreover, they are both smooth solutions of the following Stokes
problem in a half space:
\begin{equation*}
\left\{
\begin{aligned}
& -\Delta U + \na p  = 0, \quad y \in \R^2_+, \\
& \div U = 0, \quad y \in \R^2_+, \\
& U\vert_{y_2 = 0} = U_l\vert_{y_2 = 0}.
\end{aligned}
\right.
\end{equation*}
For more on the double Stokes layer potential, we refer to
\cite{Galdi:1994}. A formal energy estimate on $U_l - \tilde{U_l}$ yields,
using stationarity:
$$ \E \left( \int_{(-R,R) \times \R^+} \left|\na 
\left( U_l - \tilde{U}_l\right)(\cdot, y) \right|^2  \, dy \right) = 0, \quad \forall
R > 0. $$
One checks easily  that $\na U_l$ and $\na \tilde{U}_l$ have appropriate 
integrability, so that this estimate is
rigorous. We deduce that $U_l = \tilde{U}_l$ almost surely.

\bigskip
We now introduce 
\begin{equation} \label{Uinfty}
  U^\infty_l(\omega) = \lim_{R \rightarrow +\infty} \frac{1}{R} \int_0^R
U_l(\omega,-y_1',0) \, dy_1'. 
\end{equation}
This limit exists almost surely by the ergodic theorem, and the measurable
function $U^\infty_l$ satisfies $U_l^\infty \circ \tau_{y_1} = U_l^\infty$
for all $y_1 \in \R$. More precisely, one has   
$$U^\infty_l(\omega) = \lim_{R \rightarrow \infty} \frac{1}{R} \int_0^R
U_l(\omega,y_1-y'_1,0) \, dy_1' $$
 uniformly locally in $y_1$. Indeed, for all $|y_1| < M$, $R$ large enough,
\begin{equation*}
\begin{aligned} 
& \left| \frac{1}{R} \int_0^R
\left(U_l(\omega,y_1-y'_1,0) - U_l(\omega,-y'_1,0) \right) \, dy_1' \right|\\
& \le  \frac{1}{R} \int_{- M}^{M} | U_l(\omega,-y'_1,0 | \, dy_1' +  
\frac{1}{R} \int_{R-M}^{R+M} | U_l(\omega,-y'_1,0 | \, dy_1' \\
& \le  C(M) \frac{1}{2R}   \left( 
\int_{(-2R,2R)} | U_l(\omega,-y'_1,0 |^2 \, dy_1' \right)^{1/2} \\
& \le  C(M) \,
\frac{1}{\sqrt{R}} \left( \sup_{R\ge 1} \frac{1}{R} \int_{(-R,R)}
 | U_l(\omega,-y'_1,0 |^2 \, dy_1' \right)^{1/2} 
\end{aligned}
\end{equation*}
this last quantity 
vanishing uniformly as $R \rightarrow+\infty$, almost surely. 
Then, to  show that $U^\infty_{l,2} = 0$, we integrate $\div U_l = 0$ over 
${\cal R}^-_l(\omega,R) \cap \{y_1<0\}$:
$$0 = \int_{{\cal R}^-_l(\omega,R) \cap \{y_1<0\}} \hspace{-0.5cm} 
     \div U_l = \int_{\pa ({\cal R}^-_l(\omega,R) \cap \{y_1<0\})} 
     \hspace{-0.5cm} U_l \cdot n,$$ 
which gives, as $R \rightarrow +\infty$, 
$$   U^\infty_{l,2} = \lim_{R \rightarrow +\infty} 
\frac{1}{R} \left( \int_{\sigma(-R) \cap \{y_2 < 0 \}} \hspace{-0.5cm} 
U_{l,1} \,dy_2 - \int_{\sigma(0) \cap \{y_2 < 0 \}} \hspace{-0.5cm}
 U_{l,1} \,dy_2 \right). $$
If $U_{l,2}^\infty$ were not zero, one would  get from the left hand side
that 
$$  \int_{{\cal R}(\omega, R)} |\na U_l|^2 \ge   
C \int_{{\cal R}^-(\omega, R)} |U_l|^2 \: \ge \:  C \int_{R_0}^{R} \, dr \, 
\left( \int_{\sigma(-r) \cap \{y_2 < 0 \}}
\hspace{-0.5cm} | U_{l,1}| + \int_{\sigma(r) \cap \{y_2 < 0
  \}} \hspace{-0.5cm} | U_{l,1}| \right)^2 \: \ge \:  \delta R^3   $$
for some $\delta > 0$ and $R \ge R_0$ large enough. This would be in
contradiction with \eqref{boundsBL}b, so that $U^\infty_{l,2} = 0$. 

We now show that almost surely,
$$U_l^\infty(\omega)= \lim_{R\to+\infty} \frac{1}{R} \int_0^R 
    U_l(\omega,y'_1,0) \,dy'_1.$$
Let $\tilde U_l^\infty$ denote this limit, which exists almost surely for
the same reason as $U_l^\infty$. 
The convergence also holds in $L^2(P)$ because 
$\int_0^1 U_l(\omega,y'_1,0) \,dy'_1 \in L^2(P)$:
\begin{gather*}
\lim_{R\to +\infty} \E \left[ \left| \frac{1}{R} \int_0^R U_l(\omega,y_1,0) 
    \,dy_1 - \tilde U_l^\infty(\omega) \right|^2 \right] =0, \\
\lim_{R\to +\infty} \E \left[ \left| \frac{1}{R} \int_{-R}^0 U_l(\omega,y_1,0)
    \,dy_1 - U_l^\infty(\omega) \right|^2 \right] =0.
\end{gather*}
In particular, we have
\begin{align*}
\E[|\tilde U_l^\infty(\omega)|^2] & = \lim_{R\to +\infty} \E \left[ \left| 
   \frac{1}{R} \int_0^R U_l(\omega,y_1,0) \,dy_1 \right|^2 \right] \\
& = \lim_{R\to +\infty} \E \left[ \left| 
   \frac{1}{R} \int_{-R}^0 U_l(\omega,y_1,0) \,dy_1 \right|^2 \right] = 
   \E[|U_l^\infty(\omega)|^2]
\end{align*}
since the law of $U_l$ is invariant by translation in the $y_1$-direction.
Now, we write
\begin{align*}
\E[|\tilde U_l^\infty(\omega)|^2] & = \lim_{R\to +\infty} \E \left[ \left| 
   \frac{1}{2R} \int_0^{2R} U_l(\omega,y_1,0) \,dy_1 \right|^2 \right] \\
& = \lim_{R\to +\infty} \left( \E \left[ \left| 
   \frac{1}{2R} \int_0^R U_l(\omega,y_1,0) \,dy_1 \right|^2 \right]
   + \E \Biggl[ \left| \frac{1}{2R} \int_{R}^{2R} U_l(\omega,y_1,0) 
   \,dy_1 \right|^2 \right] \\
& \quad {} + \E \left[ \frac{1}{2R^2} \int_0^R U_l(\omega,y_1,0) \,dy_1 
   \int_{R}^{2R} U_l(\omega,y_1,0) \,dy_1 \right] \Biggr).
\end{align*}
In this expression, the first two terms tend to the same limit $\frac{1}{4}
\E[|\tilde U_l^\infty(\omega)|^2]$. Hence we deduce that the third term
converges to $\frac{1}{2} \E[|\tilde U_l^\infty(\omega)|^2]$.
Then we have 
\begin{align*}
 \E[|\tilde U_l^\infty- U_l^\infty|^2]  & = \E[|\tilde U_l^\infty|^2] +
    \E[|U_l^\infty|^2] \\
& \quad - 2 \lim_{R\to+\infty} \E \left[\frac{1}{R^2} 
    \int_0^R U_l(\omega,y_1,0) \,dy_1 \int_{-R}^0 U_l(\omega,y_1,0) \,dy_1 
    \right] \\
&  = 2 \E[|\tilde U_l^\infty|^2] -2 \lim_{R\to+\infty} \E \left[\frac{1}{R^2} 
    \int_R^{2R} U_l(\omega,y_1,0) \,dy_1 \int_0^R U_l(\omega,y_1,0) \,dy_1 
    \right] = 0.
\end{align*}

To obtain the convergence of $U_l$ to $U^\infty_l$ as $y_2 \rightarrow
+\infty$, we write by lemma \ref{doublelayer} 
\begin{equation*} 
\begin{aligned}
& U_l(\omega,y_1, y_2) - U^\infty_l(\omega) 
  = 
\int_{\R} G(y'_1,y_2) \, \left( U_l(\omega,y_1 - y_1',0) -
  U^\infty_l(\omega) \right)\,
dy'_1 \\
& = \int_{\R} G(y'_1,y_2) \, \frac{\pa}{\pa y'_1} \int_0^{y_1'} 
  U_l(\omega, y_1 - z ,0) - 
  U^\infty_l(\omega) dz \, dy_1' \\
& = -\int_{\R} \pa_1 G(y'_1,y_2) \, y_1' \, \left( \frac{1}{y'_1} 
  \int_0^{y_1'}  U_l(\omega, y_1 - z ,0) \, dz  -   U^\infty_l(\omega)
  \right) dy'_1
\end{aligned}
\end{equation*}
Let $\eps, M > 0$. There exists $R = R(M) > 0$ such that 
$$ \forall |y_1| < M, \, \forall |y'_1| > R,   
\quad  \left| \frac{1}{y'_1} 
  \int_0^{y_1'}  U_l(\omega, y_1 - z ,0) \, dz  -   U^\infty_l(\omega)
  \right| \; \le \: \eps. $$
We deduce that 
$$ \left| \int_{|y'_1| > R} \pa_1 G(y'_1,y_2) \, y_1' \, \left( \frac{1}{y'_1} 
  \int_0^{y_1'}  U_l(\omega, y_1 - z ,0) \, dz  -   U^\infty_l(\omega)
  \right) dy'_1 \right| \: \le \;  C \, \eps $$ 
Moreover, 
\begin{align*}
&  \left| \int_{|y'_1| < R} \pa_1 G(y'_1,y_2) \, y_1' \, \left( \frac{1}{y'_1} 
  \int_0^{y_1'}  U_l(\omega, y_1 - z ,0) \, dz  -   U^\infty_l(\omega)
  \right) dy'_1 \right| \\
&  \le \;  C \, 
 \left(\int_{-R-M}^{R+M}
|U_l(\omega, z,0)| dz + R |U^\infty_l(\omega)|  \right) 
 \int_{|y'_1| < R} \left| \pa_1 G(y'_1,y_2) \right| dy'_1 \\
&  \le \: C(R,M)
 \frac{1}{y_2^2} \: \xrightarrow[y_2 \rightarrow +\infty]{} 0.
\end{align*} 
This proves the almost sure convergence result. The quadratic convergence 
is  simpler and left to the reader.

%%%%%%%%%%%%%%%%%%%%%%%%%%%%%%%%%%%%%%%%%%%%%%%%%%%%%%%%%%%%%%%%%%%%%%%%
\subsection{Some more estimates}

We will also need in the sequel the following estimates.

%%%%%%%%%%%%%%%%%%%
\begin{proposition} \label{estim_deriv_BL}
We have for all $y_2>0$ and all $\alpha\in\N^2$,
\begin{equation} \label{deriv_BL1}
\sup_{R>1} \frac{1}{R} \int_{-R}^R |\pa_y^\alpha U_l(\omega,y_1,y_2)|^2 
\,dy_1 <\infty
\quad a.s.
\end{equation}
and for $|\alpha|\ge 1$,
\begin{equation} \label{deriv_BL2}
\sup_{R>1} \frac{y_2^{2|\alpha|}}{R} \E\int_{-R}^R |\pa^\alpha_y 
U_l(\cdot,y_1,y_2)|^2 
\,dy_1 \xrightarrow[y_2\to +\infty]{} 0.
\end{equation}
\end{proposition}
Note that we could also prove an almost sure version of \eqref{deriv_BL2},
as in the previous proof, but we do not need it. For $\alpha=0$, 
\eqref{deriv_BL2} holds if one replaces $U_l$ by $U_l-U_l^\infty$, as we
have seen before.\\
\textit{Proof.} 

We have
$$\pa^\alpha U_l(\omega,y_1,y_2) = \int \pa^\alpha G(y_1-y'_1,y_2) 
U_l(\omega,y'_1,0) 
\,dy'_1$$ 
hence we have $\E[|\pa^\alpha U_l(\cdot,0,y_2)|]<\infty$. The law of $\pa^\alpha 
U_l(\cdot,y_1,y_2)$ is independent of $y_1$. Thus we can apply the ergodic
theorem, 
which yields \eqref{deriv_BL1}.
We now prove \eqref{deriv_BL2}. 
We have
\begin{align*}
\pa^\alpha U_l(\omega,y_1,y_2) & = \int \pa^\alpha G(y'_1,y_2) (U_l(\omega,
y_1-y'_1,0)
     -U_l^\infty(\omega)) \,dy'_1 \\ 
& = \int y'_1 \pa_1\pa^\alpha G(y'_1,y_2) \left(\frac{1}{y'_1} \int_0^{y'_1} 
     (U_l(\omega,y_1-z,0)-U_l^\infty(\omega)) \,dz \right) \,dy'_1. 
\end{align*}
Let $M>0$. We cut the previous integral in two pieces: $|y'_1|<M$, $|y'_1|>M$,
and apply the Cauchy-Schwarz inequality in each piece. Recalling that $G$ is 
homogeneous of degree $-1$, we obtain:
\begin{align*}
|\pa^\alpha U_l(y_1,y_2)|^2 & \leq \frac{C}{y_2^{1+|\alpha|}} \int_{|y'_1|<M}
     |\pa_1\pa^\alpha G(y'_1,y_2)| \biggl| \int_0^{y'_1} (U_l(\omega,y_1-z,0)-
     U_l^\infty(\omega)) \,dz \biggr|^2 \,dy'_1 \\
& \quad {} + \frac{C}{y_2^{|\alpha|}} \int_{|y'_1|>M} |y'_1 \pa_1\pa^\alpha 
     G(y'_1,y_2)| \biggl| \frac{1}{y'_1} \int_0^{y'_1} (U_l(\omega,y_1-z,0)-
     U_l^\infty(\omega)) \,dz \biggr|^2 \,dy'_1.
\end{align*}
We now take the expectation, and obtain, due to the invariance of the
probability measure $P$ with respect to $\tau_{y_1}$,
\begin{align*}
\E[|\pa^\alpha U_l(\cdot,y_1,y_2)|^2] & \leq \frac{C}{y_2^{1+|\alpha|}} 
     \int_{|y'_1|<M} |\pa_1\pa^\alpha G(y'_1,y_2)|\, \E \biggl| \int_0^{y'_1} 
     (U_l(\cdot,-z,0)- U_l^\infty) \,dz \biggr|^2 \,dy'_1 \\
& \quad {} + \frac{C}{y_2^{|\alpha|}} \int_{|y'_1|>M} |y'_1 \pa_1\pa^\alpha 
     G(y'_1,y_2)|\, \E \biggl| \frac{1}{y'_1} \int_0^{y'_1} (U_l(\cdot,-z,0)-
     U_l^\infty) \,dz \biggr|^2 \,dy'_1.
\end{align*}
Hence the quantity in the l.h.s. is independent of $y_1$. Thus we have
\begin{align*}
\sup_{R} \frac{1}{R} \E \int_{-R}^R |\pa^\alpha U_l(\omega,y_1,y_2)|^2 \,dy_1 
     & \leq \frac{C}{y_2^{2+2|\alpha|}} \E \biggl| \int_{-M}^M (U_l(\omega,-z,0)- 
     U_l^\infty) \,dz \biggr|^2 \\
& \quad {} + \frac{C}{y_2^{2|\alpha|}}  \sup_{|y'_1|>M} \E \biggl| 
     \frac{1}{y'_1} \int_0^{y'_1} (U_l(\omega,-z,0)- U_l^\infty) \,dz 
\biggr|^2 \\
& \leq \frac{C}{y_2^{2+2|\alpha|}} + \frac{C\delta}{y_2^{2|\alpha|}}
\end{align*}
provided $M$ is chosen great enough so that $\E \bigl|\frac{1}{y'_1} 
\int_0^{y'_1} (U_l(\cdot,-z,0)- U_l^\infty) \,dz \bigr|^2<\delta$ for all 
$|y'_1|>M$. The result follows.

%%%%%%%%%%%%%%%%%%%%%%%%%%%%%%%%%%%%%%%%%%%%%%%%%%%%%%%%%%%%%%%%%%%%%%%%%%%%
%%%%%%%%%%%%%%%%%%%%%%%%%%%%%%%%%%%%%%%%%%%%%%%%%%%%%%%%%%%%%%%%%%%%%%%%%%%%
\section{Justification of Navier's law} \label{sectionNavier}

This section is divided into two parts. First we derive an approximation
$\tilde{u}^\eps_{app}$ of $u^\eps$ based on the boundary layer analysis, and prove
theorem \ref{thNavier}. Then using this approximation, we justify Navier's
wall law and prove theorem \ref{justifNavier}.

%%%%%%%%%%%%%%%%%%%%%%%%%%%%%%%%%%%%%%%%%%%%%%%%%%%%%%%%%%%%%%%%%%%%%%%%%%%
\subsection{Approximation of $u^\eps$} \label{sectionapprox}

The approximation $\tilde{u}^\eps_{app}$ reads 
\begin{equation} \label{Ansatz3}  
\begin{aligned}
\tilde{u}^\eps_{app}(\omega, x) & =  u^0(\omega, x)  + \eps  U_l\left(\omega, 
\frac{x_1}{\eps}, \frac{x_2}{\eps} \right) + \eps 
 U_u\left(\omega, \frac{x_1}{\eps},
\frac{x_2 - 1}{\eps} \right)  \\
&\quad + \eps  u^{1,\eps}(\omega,x),   
\end{aligned}  
\end{equation}
It takes into account the boundary layer
terms, and includes corrections to recover the correct boundary and flux 
conditions : this is the role of  $u^{1,\eps}$. 

More precisely, since $U_{l,u}$ does not converges to 0 when $y_2$ tends to
$\pm\infty$, we subtract $U_{l,u}^\infty$ from $U_{l,u}$ in \eqref{Ansatz}. 
Then we deal with the boundary condition. The term $\eps (U_l-U_l^\infty)$
is non zero on the upper boundary of $\Omega^\eps$. We truncate this term
replacing it by $\eps (U_l-U_l^\infty) 1_{x_2<1}$. Hence we have the
correct boundary condition, but a jump appears on $\Sigma_1$. So we
introduce a function $\eps v_l(x)$ such that
$v_l(x_1,1)=U_l^\infty-U_l(\frac{x_1}{\eps},\frac{1}{\eps})$ on 
$\Sigma_1$, $v_l=0$ on $\Sigma_0$, with $\div v_l=0$ in $\Omega$, and we set
$v_l=0$ outside $\Omega$. Notice that $\eps (U_l-U_l^\infty)=o(\eps)$, thus
$\eps v_l=o(\eps)$ (roughly), so this term will not interfere in our
estimate. 

Now on the lower boundary, we have $\eps (U_l-U_l^\infty)=-\eps
U_l^\infty$. So we will again introduce a counterflow 
$\eps c_l$ in $\Omega$. Unlike $\eps v_l$,
$\eps c_l$ is not $o(\eps)$, thus we have to choose for $c_l$ a solution 
of the following equations in $\Omega$ :
\begin{equation*} 
\left\{ \begin{aligned}
& -\nu \Delta c_l + (u^0\cdot\na)c_l + (c_l\cdot\na)u^0 + \na g_l =0 \\
& \div c_l=0 \\
& {c_l}|_{\Sigma_0}=U_{l,1}^\infty e_1, \quad {c_l}|_{\Sigma_1}=0\\
& \int_0^1 c_l \cdot e_1 \,dx_2 = 0
\end{aligned} \right.
\end{equation*}
(recall that $U_{l,2}^\infty=0$). The solution of this system is just a
combination of a Couette flow and a Poiseuille flow:
$c_l=(1-4x_2+3x_2^2) U_{l,1}^\infty e_1$, $\nabla g_l=6\nu U_l^\infty$. We
extend $c_l$ outside 
$\Omega$ by setting $c_l=0$ above $\Sigma_1$ and $c_l=U_l^\infty$ below
$\Sigma_0$. 

We proceed similarly with $\eps (U_u-U_u^\infty)$, introducing $v_u$ and
$c_u$. We have $c_u=(3x_2^2-2x_2) U_u^\infty$ in
$\Omega$, $c_u=U_u^\infty$ above $\Sigma_1$ and $c_u=0$ below $\Sigma_0$. 
Finally, we add a small Poiseuille flow $\eps \theta u_0$ to correct
the value of the flux of $u^\eps_{app}$. We obtain:
\begin{equation} \label{Ansatz4}
\tilde{u}^\eps_{app}= u^0 + \eps \Bigl( U_l\Bigl(\frac{x}{\eps}\Bigr)
-U_l^\infty \Bigr) 1_{x_2<1} + \eps \Bigl( U_u\Bigl(\frac{x}{\eps}\Bigr)
-U_u^\infty \Bigr) 1_{x_2>0} +\eps(v_l+v_u+c_l+c_u) + \eps\theta u^0
\end{equation}
where the coefficient $\theta$ is defined by
$$\theta \phi = -\int_{\sigma(x_1)} \Bigl(
U_{l,1}\Bigl(\frac{x}{\eps}\Bigr) 
1_{x_2<1} + U_{u,1}\Bigl(\frac{x}{\eps}\Bigr) 1_{x_2>0} \Bigr) \,dx_2 +
U^\infty_{l,1} + U^\infty_{u,1} - \int_0^1 (v_{l,1}(x) + v_{u,1}(x)) \,dx_2.$$ 
The value of $\theta$ is independent of $x_1$, because $\div
\tilde u^\eps_{app}=0$. Easy recombinations show that the expansion
\eqref{Ansatz4} is of type \eqref{Ansatz3}, with 
$$u^{1, \eps} \:  = \:  u^1 \: + \: v_l \: +\: v_u \: + \: \theta \, u^0
 - \Bigl( U_l\Bigl(\frac{x}{\eps}\Bigr)
-U_l^\infty \Bigr) 1_{x_2>1} + \eps \Bigl( U_u\Bigl(\frac{x}{\eps}\Bigr)
-U_u^\infty \Bigr) 1_{x_2<0},$$
and $u^1=c_l+c_u-U_l^\infty-U_u^\infty$ (according to \eqref{Ansatz2} and
\eqref{u1}). 

%%%%%%%%%%%%%%%%%%%%%%%%%%%%%%%%%%%%%%%%%%%%%%%%%%%%%%%%%%%%%%
\subsubsection{Construction of $v_l$ and $v_u$}

We explain in this paragraph the construction of the flow $v_l$. The
construction of $v_u$ is analogous.
We have to solve the following problem: find $v_l \in H^1_{loc}(\Omega)$
such that
\begin{equation*}
\left\{\begin{aligned}
& v_l(x_1,0)=0 \text{ on } \Sigma_0 \\
& v_l(x_1,1)=U_l^\infty - U_l(\tfrac{x_1}{\eps},\tfrac{1}{\eps}) \text{ on
       } \Sigma_1 \\
& \div v_l=0 \text{ in } \Omega.
\end{aligned}\right.
\end{equation*}

\begin{proposition} \label{exist_vl}
This problem possesses a (non unique) solution $v_l$ such that
$$\sup_{R>1} \frac{1}{R} \|v_l\|_{H^2(\Omega_R)}^2 < \infty \quad a.s.$$
and %%%for all $R>1$,
$$\sup_{R>1} \frac{1}{R} \E \bigl(\|v_l\|_{H^2(\Omega_R)}^2 \bigr) = o(1)$$
when $\eps \to 0$.
\end{proposition}
\textit{Proof.} 

We shall find a solution in the following form: $v_l=\na^\perp \psi$. In terms
of $\psi$, the boundary conditions can be rewritten as $\pa_1 \psi=\pa_2 \psi
=0$ on $\Sigma_0$, $\pa_1 \psi(x_1,1)=-U_{l,2}(\frac{x_1}{\eps},
\frac{1}{\eps})$ and $\pa_2 \psi(x_1,1)=U_{l,1}(\frac{x_1}{\eps},
\frac{1}{\eps}) - U_{l,1}^\infty$ on $\Sigma_1$. Up to a constant, this is
equivalent to $\psi=\pa_2\psi=0$ on $\Sigma_0$, $\psi(x_1,1)=-\int_0^{x_1}
U_{l,2}(\frac{x'_1}{\eps},\frac{1}{\eps}) \,dx'_1$, $\pa_2 \psi(x_1,1)=
U_{l,1}(\frac{x_1}{\eps},\frac{1}{\eps}) - U_{l,1}^\infty$ on $\Sigma_1$.
We search a solution which is polynomial in $x_2$:
let $\psi(x_1,x_2)=a(x_1)x_2^3+b(x_1)x_2^2+c(x_1)x_2+d(x_1)$. The boundary
conditions on $\Sigma_0$ imply $c(x_1)=d(x_1)=0$, and on $\Sigma_1$ we get
\begin{equation*}
\left\{\begin{aligned}
a+b & = -\int_0^{x_1} U_{l,2}\Bigl(\frac{x'_1}{\eps},\frac{1}{\eps}\Bigr) 
        \,dx'_1 \\
3a+2b & = U_{l,1}\bigl(\frac{x_1}{\eps},\frac{1}{\eps}\bigr) - U_{l,1}^\infty.
\end{aligned}\right.
\end{equation*}
This system has of course a unique solution $(a(x_1),b(x_1))$, and to obtain
the required estimates, it is sufficient to show that
\begin{gather*}
\sup_{R>1} \frac{1}{R}\left(\Bigl\|\int_0^{x_1} U_{l,2}(\tfrac{x'_1}{\eps},
        \tfrac{1}{\eps}) \,dx'_1 \Bigr\|_{H^3(-R,R)}^2 + \bigl\|
        U_{l,1}(\tfrac{x_1}{\eps},\tfrac{1}{\eps}\bigr) - U_{l,1}^\infty
        \bigr\|_{H^3(-R,R)}^2 \right) < \infty \quad \text{a.s.}, \\
\sup_{R>1} \frac{1}{R} \E \left[\Bigl\|\int_0^{x_1} 
        U_{l,2}(\tfrac{x'_1}{\eps},\tfrac{1}{\eps}) \,dx'_1 
        \Bigr\|_{H^3(-R,R)}^2 + \bigl\|U_{l,1}(\tfrac{x_1}{\eps},
        \tfrac{1}{\eps}\bigr) - U_{l,1}^\infty\bigr\|_{H^3(-R,R)}^2 \right] 
        \to 0.
\end{gather*}
The estimate on the quantity $\|U_l(\frac{x_1}{\eps},
\frac{1}{\eps}\bigr) - U_l^\infty\|_{H^3(-R,R)}^2$ follows from  proposition 
\ref{estim_deriv_BL}, for we have
\begin{align*}
\frac{1}{R} \bigl\|U_l(\tfrac{x_1}{\eps},\tfrac{1}{\eps}\bigr) - U_l^\infty
      \bigr\|_{H^3(-R,R)}^2 = \frac{\eps}{R} \int_{-R/\eps}^{R/\eps}
      \bigl( & |U_l(y_1,\tfrac{1}{\eps})-U_l^\infty|^2 + \eps^{-2} |\pa_1
      U_l(y_1,\tfrac{1}{\eps})|^2\\
& {} + \eps^{-4} |\pa_1^2  
      U_l(y_1,\tfrac{1}{\eps})|^2  +  \eps^{-6} |\pa_1^3  
      U_l(y_1,\tfrac{1}{\eps})|^2\bigr) \,dy_1
\end{align*}
(recall that $U_{l,2}^\infty=0$).
Hence we only need to deal with $\frac{1}{R} \bigl\| \int_0^{x_1} 
U_{l,2}(\frac{x'_1}{\eps},\frac{1}{\eps}) \,dx'_1 \bigr\|_{L^2(-R,R)}^2$.
In order to do so, we write:
\begin{align*}
\int_0^{x_1} U_{l,2}(\tfrac{x'_1}{\eps},\tfrac{1}{\eps}) \,dx'_1 & =
    \eps \int_0^{x_1/\eps} U_{l,2}(y_1,\tfrac{1}{\eps}) \,dy_1 \\
& = \eps \int_{\sigma(-\frac{x_1}{\eps}) \cap \{y_2<\frac{1}{\eps}\}}
    \hspace{-0.5cm} U_{l,1}(\tfrac{x_1}{\eps},y_2) \,dy_2 - \eps 
    \int_{\sigma(0) \cap \{y_2<\frac{1}{\eps}\}} \hspace{-0.5cm}
    U_{l,1}(0,y_2) \,dy_2
\end{align*}
since $U_l$ is divergence free on $\mathcal{R}_l \cap \{0<y_1<
\frac{x_1}{\eps},\ y_2<\frac{1}{\eps}\}$. 
Now the pointwise estimate is easy, for we have
$$\E\biggl| \eps \int_{\sigma(0) \cap \{y_2<\frac{1}{\eps}\}}
    \hspace{-0.5cm} U_{l,1}(0,y_2) \,dy_2 \biggr|^2 < \infty,$$
hence we can apply the pointwise ergodic theorem. 
Next we turn to the estimate in expectation. 
We have
\begin{align*}
\int_0^{x_1} U_{l,2}(\tfrac{x'_1}{\eps},\tfrac{1}{\eps}) \,dx'_1 & =
    \eps\int_0^{1/\eps} (U_{l,1}(\tfrac{x_1}{\eps},y_2)-U_{l,1}^\infty) \,dy_2
    - \eps \int_0^{1/\eps} (U_{l,1}(0,y_2)-U_{l,1}^\infty) \,dy_2 \\
& \quad {} + \eps \int_{\sigma(-\frac{x_1}{\eps}) \cap \{y_2<0\}}
    \hspace{-0.5cm} U_{l,1}(\tfrac{x_1}{\eps},y_2) \,dy_2 - \eps 
    \int_{\sigma(0) \cap \{y_2<0\}} \hspace{-0.5cm}
    U_{l,1}(0,y_2) \,dy_2
\end{align*}
Then the invariance of $P$ yields 
$$
\E \frac{1}{2R} \biggl\| \eps\int_0^{1/\eps} (U_{l,1}(\tfrac{x_1}{\eps},y_2)-
    U_{l,1}^\infty) \,dy_2 \biggr\|_{L^2(-R,R)}^2 = \E \biggl| \eps
    \int_0^{1/\eps} (U_{l,1}(0,y_2)- U_{l,1}^\infty) \,dy_2 \biggr|^2.$$
Therefore it is sufficient to deal with the integral of the l.h.s. with $R=1$.
Let $\delta>0$. We have
\begin{align*}
\E \biggl\| \eps\int_0^{1/\eps} (U_{l,1}(\tfrac{x_1}{\eps},y_2)-
    U_{l,1}^\infty) \,dy_2 \biggr\|_{L^2(-1,1)}^2 & \leq \eps\int_0^{1/\eps} 
    \E \biggl[\int_{-1}^1 |U_{l,1}(\tfrac{x_1}{\eps},y_2)-
    U_{l,1}^\infty|^2 \,dx_1\biggr] \,dy_2 \\
& \leq \eps\int_0^{1/\eps} \E \biggl[\int_{-1}^1 |U_{l,1}(y_1,y_2)-
    U_{l,1}^\infty|^2 \,dy_1\biggr] \,dy_2 \\
& \leq \eps\E \int_{[{-1},1]\times [0,M]} |U_{l,1}-U_{l,1}^\infty|^2 + \eps 
    \int_M^{1/\eps} \delta \,dy_2 \\
& \leq C\eps +\delta,
\end{align*}
where we have used successively the Cauchy-Schwarz inequality, the invariance 
of $P$ and the estimate \eqref{cvgce_BL} to make the expectation of 
$\int_{-1}^1 |U_{l,1}-U_{l,1}^\infty|^2 \,dy_1$ less than $\delta$ if $y_2$ 
is larger than some constant $M$.
Letting $\eps$ decrease to 0, we deduce that the l.h.s. above tends to 0.
It remains to deal with the integrals on $\sigma(\frac{x_1}{\eps}) \cap
\{y_2<0\}$ and $\sigma(0) \cap \{y_2<0\}$. As above, the expectations of
the quadratic mean on $(-R,R)$ of these two integrals are equal, and equal to:
\begin{align*}
\E \frac{1}{2R} \biggl\| \eps \int_{\sigma(-\frac{x_1}{\eps}) \cap \{y_2<0\}}
    \hspace{-0.5cm} U_{l,1}(\tfrac{x_1}{\eps},y_2) \,dy_2 
    \biggr\|_{L^2(-R,R)}^2 & = \E \biggl|\eps \int_{\sigma(0) \cap \{y_2<0\}} 
    \hspace{-0.5cm} U_{l,1}(0,y_2) \,dy_2\biggr|^2 \\
& = \E \int_0^1 \biggl|\eps \int_{\sigma(y_1) \cap \{y_2<0\}} \hspace{-0.5cm}
    U_{l,1}(y_1,y_2) \,dy_2\biggr|^2 \,dy_1 \\
& \leq \eps^2 \E \int_0^1 \int_{\sigma(y_1) \cap \{y_2<0\}} \hspace{-0.5cm}
    |U_{l,1}(y_1,y_2)|^2 \,dy_2 dy_1 \\
& \leq C\eps^2.
\end{align*}
Putting together the four previous estimates, we can conclude that
$$\sup_R \frac{1}{R} \E \biggl\| \int_0^{x_1} U_{l,2}(\tfrac{x'_1}{\eps},
\tfrac{1}{\eps}) \,dx'_1 \biggr\|_{L^2(-R,R)}^2 \to 0.$$
This ends the proof of proposition \ref{exist_vl}.

\subsubsection{Proof of theorem \ref{thNavier}}
This paragraph is devoted to the estimates \eqref{estimatesNavier}. It is
enough to derive the first inequality, as the second one follows 
 directly from Poincar\'e inequality. In fact, we will  show that 
$$  \sup_{R \ge 1} \frac{1}{R} \, \E \left( \int_{\Omega^\eps(\cdot,R)}
 \left|\na u^\eps -
 \na \tilde{u}^\eps_{app}\right|^2 \right) \: = \: o(\eps^2) $$
 where $\tilde{u}^\eps_{app}$ is given by \eqref{Ansatz3}-\eqref{Ansatz4}. 
Indeed, one verifies
easily from the previous estimates that 
$$ \sup_{R \ge 1} \frac{1}{R} \, \E \left( \int_{\Omega^\eps(\cdot,R)}
 \left|\na u^\eps_{app} -
 \na \tilde{u}^\eps_{app}\right|^2 \right) \: = \: o(\eps^2). $$
We write 
$\dis  w = u^\eps - \tilde{u}^\eps_{app} = v^\eps - v^\eps_{app}, \: v^\eps =
u^\eps - u^0, \: v^\eps_{app} = \tilde{u}^\eps_{app} - u^0. $
It satisfies \eqref{linearNS} with $v = 0$, 
\begin{equation*}
\left\{
\begin{aligned}
  f(\omega, x)  & = 0, \quad x \mbox{ in } \Omega, \\
 f(\omega, x)   & = \biggl( -12\nu \phi \, (1 + \eps \theta),
0 \biggr), \quad x \mbox{ in } \Omega^\eps(\omega) \setminus \Omega, \\
\end{aligned}
\right.
\end{equation*}
and
\begin{equation*}
\tilde{\phi}\vert_{\Sigma_{0,1}}  \:  =  \:\pm \left[ \pa_2
  U_{u,l}\left(\frac{x}{\eps}\right)\right]\vert_{\Sigma_{0,1}}  \:+ \: \eps
\, \pa_2 (v_l + v_u + c_l + c_u)\vert_{\Sigma_{0,1}} \: +  \:6 \eps \theta
\phi.
\end{equation*}
Finally, using that  
$$\div \left(u^0 \otimes \left(U^\infty_l + U^\infty_u + u^1\right)\right)
 = \div \left(
\left(U^\infty_l + U^\infty_u + u^1 \right) \otimes u^0 \right) = 0,$$ 
where $u^1$ is given in \eqref{u1}, we write
\begin{align*} 
  G  \:  & = \: -v^\eps \otimes v^\eps  \:+  \:u^0 \otimes \biggl(
v^\eps_{app} - \eps \left(U^\infty_l +  U^\infty_u + u^1\right) \biggr) \\
& \quad   + \: \biggl(v^\eps_{app} - \eps \left( U^\infty_l + U^\infty_u +
u^1\right)  \biggr)   \otimes u^0   \:+ \: \nu \eps \na (v_l + v_u).
\end{align*} 
Proceeding as in section \ref{sectionDirichlet}, one has
\begin{align*}
 \int_\eta^{\eta+1} dR \, \int_{\Omega^\eps(R)} |\na w|^2 \: & \le \: 
 C \biggl(  
\int_\eta^{\eta+1} dR  \, \int_{\Omega^\eps(R)} \left( | G |^2 + \eps^2 | f |^2
\right) \: + \: \eps \, \int_\eta^{\eta+1} dR  \int_{\Sigma_0(R) \cup \Sigma_1(R)} 
| \tilde{\phi} |^2 \\ 
& \quad + \: \int_{\Omega^\eps(\eta, \eta+1)}\left( |  G |^2 + \eps^2 
| f |^2 \right) \: + \: \int_{\Omega^\eps(\eta, \eta+1)} | \na w |^2 \biggr)
\end{align*}
which yields roughly, after integration with respect to $\omega$:
\begin{align*}
\E  \, \int_\eta^{\eta+1} dR \, \left( \int_{\Omega^\eps(R)} 
|\na w|^2 \right) \: & \le \: 
 C \biggl( \E 
\int_{\Omega^\eps(\eta+1)} \hspace{-0.5cm} \left( | G |^2 + \eps^2 | f |^2
\right)  \: + \: \eps \,  \E  \int_{\Sigma_0(\eta+1) \cup
  \Sigma_1(\eta+1)}  \hspace{-0.5cm} 
| \tilde{\phi} |^2    \\
  & \quad + \:  \E \int_{\Omega^\eps(\eta, \eta+1)} | \na w |^2 \biggr).
\end{align*}
Assuming that 
\begin{equation} \label{lastestimates}
\begin{aligned}
& \E \int_{\Omega^\eps(\eta+1)}  \hspace{-0.5cm} 
|G|^2 +  \eps^2 | f |^2 \: \le  \:  \delta(\eps) \, \eps^2
  \left( \eta + 1 \right), \quad \delta(\eps) \xrightarrow[\eps
  \rightarrow 0]{} 0, 
 \\
& \E  \int_{\Sigma_{0,1}(\eta+1)} \hspace{-0.5cm}
| \tilde{\phi} |^2  \: \le \:  C \, \eps^2 \, (\eta + 1), 
\end{aligned}
\end{equation}
we end up  with a reverse Gronwall inequality:
$$   F(\eta) \: \le \: \eta(\eps) \, \eps^2 (\eta + 1) \: + \:  C \, F'(\eta), 
  \quad F(\eta) =\E  \, \int_\eta^{\eta+1} dR \, \left( \int_{\Omega^\eps(R)} 
 |\na w|^2 \right), \quad \eta(\eps) \xrightarrow[\eps
  \rightarrow 0]{} 0, $$ 
and conclude as in section \ref{sectionDirichlet}. It thus remains to
establish bounds \eqref{lastestimates}. The second inequality is obvious
using propositions \ref{estim_deriv_BL} and \ref{exist_vl}. To control $G$,
we first notice that 
\begin{align*}
\int_{\Omega^\eps(\eta+1)} | v^\eps \otimes v^\eps |^2  \: & \le \: 
 C \, (\eta +  1) \, \sup_{ |k| \le \eta + 2} \| v^\eps
 \|^{4}_{L^4(\Omega^\eps(k,k+1))} \\
& \le \:  C \, (\eta +  1) \, \| v^\eps \|_{L^2_{uloc}}^2 \, \| \na v^\eps
 \|_{L^2_{uloc}}^2 
\end{align*}
by standard Gagliardo-Nirenberg inequality.  Using estimates
\eqref{Dirichlet_uloc1}, \eqref{Dirichlet_uloc3}, we deduce 
$$ \int_{\Omega^\eps(\eta+1)} | v^\eps \otimes v^\eps |^2  \:  \le \: 
 C \, \eps^3 \, (\eta +  1). $$
Then, with notations of \eqref{Ansatz3} we write 
\begin{align*}
&  v^\eps_{app}(\omega,x) - \eps \left(U_l^\infty(\omega) -
U_u^\infty(\omega) - u^1(\omega,x)\right) \:  = \: \eps \biggl( 
\left( U_l^\infty\left(\omega,\frac{x}{\eps}\right) -
U_l^\infty(\omega)\right) \\
&   + \eps 
\left( U_u^\infty\left(\omega,\frac{x_1}{\eps}, \frac{x_2-1}{\eps}\right) -  
U_u^\infty(\omega)\right)  + \eps \left( u^{1,\eps}(\omega,x) -
u^1(\omega,x) \right)\biggr). 
\end{align*}
In particular, the estimate of $G$ involves
\begin{align*}
& \E \left(\int_{\Omega^\eps(\omega,\eta+1)} 
 \left| U_u^\infty\left(\omega,\frac{x_1}{\eps}, \frac{x_2-1}{\eps}\right) -  
U_u^\infty(\omega)\right| dx \right)  
 \, d\omega  \\  
& \le \: C \, (\eta + 1) \, 
\E \left( \eps \int_{-h_l(0)}^{1/\eps + h_u(0)} \left|
 U_l(\omega,0, y_2)  - U^l_{\infty}(\omega) \right|^2 dy_2 \right) \, d\omega =
 o(1), \quad \eps \rightarrow 0,
\end{align*}
with  the same  manipulations as in the proof of proposition
\ref{exist_vl}. 
All other terms lead to similar inequalities and are left to the reader. This
ends the proof. 
\subsection{Effective wall law}
It remains to connect the approximation we have  built to the Navier's
wall law. Following \cite{Jager:2001}, 
such connection can be seen at a formal level,  as
\begin{equation*} 
u^\eps  \approx u^\eps_{app} =  u^0(\omega, x)  + \eps  U_l\left(\omega, 
\frac{x_1}{\eps}, \frac{x_2}{\eps} \right) 
 + \eps U_u\left(\omega, \frac{x_1}{\eps},
\frac{x_2 - 1}{\eps} \right)  \:  + \eps  u^1(\omega,x),
\end{equation*}
For instance, we obtain formally: 
$$ u^\eps_1\vert_{\Sigma_0} \: \approx \: \eps U_{l,1} \left(\omega,
\frac{x_1}{\eps}, 0\right), \quad \pa_2 u^\eps_1\vert_{\Sigma_0} \: \approx
\: 6 \phi + \pa_2 U_{l,1}\left(\omega,\frac{x_1}{\eps}, 0\right), $$
which yields after averaging the approximate boundary condition,
$$ 
 v^\eps_1\vert_{\Sigma_0} =  \eps \alpha_l  \, 
\pa_2 v^\eps_1\vert_{\Sigma_0}. 
$$
Similarly, 
$$
 v^\eps_1\vert_{\Sigma_1} =  -\eps \alpha_u \, \pa_2 
v^\eps_1\vert_{\Sigma_1}. 
$$

We now justify this formal computation and prove theorem
\ref{justifNavier}. 
 Let $v^\eps$ the solution of \eqref{Naviereq}-\eqref{Navierlaw}. It
is given explicitly by  
$$ v^\eps = \biggl(M_\phi \, \left( x_2^2 - (x_2 - \eps \alpha_l)
\frac{1-2\eps \alpha_u}{2(1-\eps (\alpha_l + \alpha_u))}\right), \, 0
\biggr), \quad \phi \, M_\phi 
 = \frac{1}{3} - \frac{(1-2\eps \alpha_l)(1 -2 \eps 
  \alpha_u)}{1-\eps(\alpha_l + \alpha_u)}.$$      
Then, we use the identity:   
\begin{align*}
u^\eps(\omega,x) - v^\eps(\omega,x)  \:  & = \: 
 u^\eps(\omega,x) - u^\eps_{app}(\omega,x) + u^\eps_{app}(\omega,x) -
 v^\eps(\omega,x)  =  \left(u^\eps - u^\eps_{app}\right)(\omega,x) \\
& \quad  + \eps \biggl( u^0(\omega,x) + \eps u^1(\omega,x)  -
 U^\infty_l(\omega)  - U^\infty_u(\omega) \biggr) 
-  v^{\eps}(\omega,x)  \\
& \quad + \eps  \left( U_l\left(\omega, 
\frac{x_1}{\eps}, \frac{x_2}{\eps} \right) - U_l^\infty(\omega)\right)
 + \eps \left( U_u\left(\omega, \frac{x_1}{\eps}, 
\frac{x_2 - 1}{\eps} \right) - U_u^\infty(\omega) \right). 
\end{align*}
By theorem \ref{thNavier}, we have 
$$ \sup_{R \ge 1} \frac{1}{R} \, \E \left( 
\int_{\Omega(R)} \left|  u^\eps(\omega, \cdot) 
- u^\eps_{app}(\omega,\cdot) \right|^2  \right)  \: =  \:  o(\eps^2). $$
Then, simple calculations show that the function
 $$w =  \eps \left( u^0(\omega,x) + \eps u^1(\omega,x)  -
 U^\infty_l(\omega)  - U^\infty_u(\omega) \right)  -  v^{\eps}(\omega,x) $$
which is the solution of 
\begin{equation*}
\left\{
\begin{aligned}
& w \cdot \na w  - \nu \Delta w + \na q = 0, \quad x \in \Omega, \\
& \div w = 0,  \quad x \in \Omega, \\
&  \int_{\sigma(x_1)} \!\!\!\! w_1 \:= \: 0,  \quad w_2(\omega,
  \cdot)\vert_{\Sigma_{0,1}} = 0, \\
& w(\omega, \cdot)   - \eps \, \alpha_l(\omega) \,    \, 
\frac{\pa w(\omega, \cdot)}{\pa x_2} = \eps^2 \alpha_l(\omega) \left(4
U^\infty_{l,1}(\omega) + 2 U^\infty_{u,1}(\omega)\right),   \quad x_2 = 0, \\
& w(\omega, \cdot) + \eps \alpha_u(\omega)  \, \frac{\pa
    w(\omega, \cdot)}{\pa x_2} = \eps^2 \alpha_u(\omega)\left(2
U^\infty_{l,1}(\omega) + 4 U^\infty_{u,1}(\omega)\right) ,  \quad x_2 = 1, 
\end{aligned}
\right.
\end{equation*}
satisfies 
$$ \sup_{R \ge 1} \frac{1}{R} \, \E \left(  
\int_{\Omega(R)} \left|  u^\eps(\omega, \cdot) 
- u^\eps_{app}(\omega,\cdot) \right|^2  \right)  \: =  \:  O(\eps^4). $$
Finally, we proceed as usual to get 
$$ \sup_{R \ge 1} \frac{1}{R} \, \E \left(  
\int_{\Omega(R)}  \left|   U_l\left(\omega, 
\frac{x_1}{\eps}, \frac{x_2}{\eps} \right) - U_l^\infty(\omega) \right|^2 dx
\right) d\omega  = o(1), $$
and similar estimate for the upper boundary layer. This concludes the proof
of theorem \ref{justifNavier}. 

%%%%%%%%%%%%%%%%%%%%%%%%%%%%%%%%%%%%%%%%%%%%%%%%%%%%%%%%%%%%%%%%%%%%%%%%%%
%%%%%%%%%%%%%%%%%%%%%%%%%%%%%%%%%%%%%%%%%%%%%%%%%%%%%%%%%%%%%%%%%%%%%%%%%%

\bibliographystyle{plain} 

%%%%%%%%%%%%%%%%

\end{document}